\documentclass[letterpaper,12pt]{amsart}
\usepackage{amsfonts,color,latexsym,amssymb,amsmath,epsfig,rotating,array,vicent}
\usepackage[OT1]{fontenc} 

\newtheorem{adelic}{Adelic $\pmb{\tau}$-Conjecture}

\newtheorem{kapra}{Kapranov's Non-Archimedean Amoeba Theorem}

\newtheorem{dfn}{Definition}[section]
\newtheorem{rem}[dfn]{Remark} 
\newtheorem{prop}[dfn]{Proposition}

\newtheorem{thm}[dfn]{Theorem}
 
\newtheorem{lemma}[dfn]{Lemma}

\newtheorem{ex}[dfn]{Example}
\definecolor{dblue}{rgb}{0,0,.3}

\newcommand{\floor}[1]{\left\lfloor#1\right\rfloor} 
\newcommand{\ceil}[1]{\left\lceil#1\right\rceil}

\newtheorem{localconj}{The Local Fewnomial Conjecture}

\newtheorem{sturmf}{Sturmfels' Theorem for Complete Intersections 
(special case)}

\renewcommand{\qed}{$\blacksquare$}
\newcommand{\bg}{\bar{g}}
\newcommand{\jac}{\mathrm{Jac}}
\newcommand{\sps}{\mathrm{SPS}}

\newlength{\burg}
\newlength{\hwl}
\newlength{\koi}
\newlength{\sma}
\newlength{\jmr}
\newlength{\khov}
\newcommand{\thth}{^{\text{\underline{th}}}}

\newcommand{\stst}{^{\text{\underline{st}}}}
\newcommand{\ord}{\mathrm{ord}}

\newcommand{\np}{{\mathbf{NP}}}

\newcommand{\fk}{{\textcursive{k}}}
\newcommand{\vol}{{\mathrm{Vol}}}

\newcommand{\pp}{\mathbf{P}}
\newcommand{\vp}{\mathbf{VP}}
\newcommand{\vnp}{\mathbf{VNP}}

\newcommand{\eps}{\varepsilon}

\newcommand{\F}{\mathbb{F}}
\newcommand{\Q}{\mathbb{Q}}
\newcommand{\R}{\mathbb{R}}

\newcommand{\C}{\mathbb{C}}
\newcommand{\fM}{\mathfrak{M}}
\newcommand{\N}{\mathbb{N}}
\newcommand{\Z}{\mathbb{Z}}
\newcommand{\bO}{\mathbf{O}}

\newcommand{\cp}{\mathfrak{p}}

\newcommand{\Zn}{\Z^n}

\newcommand{\Rn}{\R^n}

\newcommand{\Cn}{\C^n}

\newcommand{\Cs}{\C^*}
\newcommand{\Rs}{\R^*}
\newcommand{\cA}{{\mathcal{A}}} 
\newcommand{\bA}{\bar{\cA}} 
\newcommand{\hA}{\hat{\cA}} 
\newcommand{\hE}{\hat{E}} 
\newcommand{\cM}{\mathcal{M}} 
\newcommand{\hP}{\hat{P}} 
\newcommand{\hQ}{\hat{Q}} 
\newcommand{\hT}{\hat{T}}

\newcommand{\bL}{{\bar{L}}}

\newcommand{\cV}{{\mathcal{V}}}

\newcommand{\Csn}{{(\C^*)}^n}

\newcommand{\dia}{$\diamond$}

\newcommand{\trop}{\mathrm{Trop}}

\newcommand{\newt}{\mathrm{Newt}}
\newcommand{\supp}{\mathrm{Supp}}
\newcommand{\conv}{\mathrm{Conv}}

\newcommand{\sign}{\mathrm{sign}}

\author{Kaitlyn Phillipson$^*$}\thanks{$^*$n\'ee Hellenbrand}
\address{Department of Mathematics,
Texas A\&M University TAMU 3368,
College Station, Texas \ 77843-3368, USA. }
\email{kaitlyn@math.tamu.edu} 
\author{J.\ Maurice Rojas}
\email{rojas@math.tamu.edu} 
\thanks{K.P.\ and J.M.R.\ were partially supported by NSF MCS 
grant DMS-0915245 and DOE ASCR grant DE-SC0002505. J.M.R.\ was also 
partially supported by Sandia National Laboratories.}  
\title[Fewnomial Systems with Many Roots, and Adelic Tau]{\mbox{}\\  
\vspace{-1in}Fewnomial Systems with Many Roots, and an Adelic Tau 
Conjecture}

\begin{document}
\maketitle 

\mbox{}\hfill
{\em To Bernd Sturmfels on his $50\thth$ birthday.}   
\hfill\mbox{}

\begin{abstract}  
Consider a system $F$ of $n$ polynomials in $n$ variables, with a 
total of \mbox{$n+k$} distinct exponent vectors, over any local field $L$.  
We discuss conjecturally tight 
bounds on the maximal number of non-degenerate
roots $F$ can have over $L$, 
with all coordinates having fixed phase, as a function of 
$n$, $k$, and $L$ only. 
In particular, we give new explicit systems with number of roots  
approaching the best known upper bounds. 
We also briefly review the background 
behind such bounds, and their application, including connections to 
computational number theory and variants 
of the Shub-Smale $\tau$-Conjecture and the $\pp$ vs.\ $\np$ Problem. 
One of our key tools is the construction of 
combinatorially constrained tropical varieties with maximally many 
intersections. 
\end{abstract} 

\section{Introduction}  
Let $L$ be any local field, i.e., $\C$, $\R$, any 
finite algebraic extension of $\Q_p$, or $\F_q((t))$. Also let 
$f_1,\ldots,f_n\!\in\!L\!\left[x^{\pm 1}_1,\ldots,
x^{\pm 1}_n\right]$ be Laurent polynomials such that the total number of 
distinct exponent vectors in the monomial term expansions of $f_1,\ldots,f_n$ 
is $n+k$. We call $F\!:=\!(f_1,\ldots,f_n)$ an {\em 
$(n+k)$-nomial $n\times n$ system over $L$}.  
We study the distribution of the 
non-degenerate roots\footnote{i.e., roots with Jacobian of rank $n$} of 
$F$ in the multiplicative group $(L^*)^n$, as a function of $n$, $k$, and $L$ 
only. This is a fundamental problem in {\em fewnomial theory over local 
fields}. We will sometimes refer to the cases $L\!\in\!\{\R,\C\}$
as the {\em Archimedean} case. Our main focus will be 
the number of roots in a fixed angular direction from the origin. 

Fewnomial theory over $\R$ has since found applications in Hilbert's 16$\thth$
Problem \cite{kaloshin}, the complexity of geometric algorithms 
\cite{gv,butterfly,thresh,prt,batessottile,bhpr,koiran,koiwron},
model completeness for certain theories of real analytic functions
\cite{wilkie,servi}, and the study of torsion points on curves 
\cite{cohenzannier}. Fewnomial theory over number fields has applications to
sharper uniform bounds on the number of torsion points on elliptic curves
\cite{cheng}, integer factorization \cite{lipton}, additive complexity
\cite{add}, and polynomial factorization and interpolation 
\cite{lenstra2,kaltofen,avenfactor,giesbrecht,newkoiran}. 
In Section \ref{sec:app} we also present an application of general fewnomial 
bounds to circuit complexity. Since any number field embeds in some finite 
extension of $\Q_p$, we thus have good reason to study fewnomial bounds over 
non-Archimedean fields. However, for $n,k\!\geq\!2$, {\em tight} bounds remain 
elusive \cite{tri,amd,bs,ai,ai2}.

\begin{dfn}
\label{dfn:omega} 
Let $y\!\in\!L^*$. When $L\!\in\!\{\R,\C\}$ we let
$|y|$ denote the usual absolute value and define
$\phi(y)\!:=\!\frac{y}{|y|}$ to be the {\em generalized phase} of
$y$. In the non-Archimedean case, we let $\fM$ denote the unique
maximal ideal of the ring of integers of $L$ and call any generator $\rho$
of $\fM$ a {\em uniformizer} for $L$. Letting $\ord$ denote the
corresponding valuation on $L$ we then alternatively define the generalized
phase as $\phi(y)\!:=\!\frac{y}{\rho^{\ord \; y}} \text{ mod } \fM$. Finally,
for general local $L$, we define $Y_L(n,k)$ 
to be the supremum, over all $(n+k)$-nomial $n\times n$ systems $F$ over $L$, 
of the number of non-degenerate
roots of $F$ in $L^n$ with all coordinates having generalized phase $1$. \dia
\end{dfn}

\noindent
Note that $y\!\in\!\C$ has generalized phase $1$ if and only if $y$ is 
positive. In the non-Archimedean case, $\phi(y)$ can be regarded simply as the 
first digit of an expansion of $y$ as a Laurent series in $\rho$. 
It is well-known in number theory that $\phi(y)$ is a natural extension of 
the argument (or angle with respect to the positive ray) of a complex 
number.\footnote{See, e.g., Schikhof's notion of {\em sign group} in 
\cite[Sec.\ 24, pp.\ 65--67]{schikhof}.}  
Our choices of uniformizer and angular direction above are in fact immaterial 
for the characteristic zero case: see Proposition \ref{prop:iou} of Section 
\ref{sec:wrap}, which also discusses the positive characteristic case. 

Descartes' classic 17$\thth$ century bound on
the number of positive roots of a sparse (a.k.a.\ lacunary) univariate 
polynomial \cite{descartes,wang}, along with some late to post-20th century 
univariate bounds of Voorhoeve, H.\ W.\ Lenstra (Jr.), Poonen, Avendano, and 
Krick, can then be recast as follows:  
\begin{thm} 
\label{thm:sum} 
Let $p$ be prime and $k\!\geq\!1$. Then:  
(1) $Y_\R(1,k)\!=\!k$ and $Y_\C(1,k)\!=\!k$, \linebreak 
(2) $Y_{\Q_2}(1,1)\!=\!2$, 
(3) $Y_{\Q_2}(1,2)\!=\!6$, 
(4) $Y_{\Q_p}(1,1)\!=\!1$ for $p\!\geq\!3$, 
(5) $Y_{\Q_p}(1,2)\!=\!3$ for $p\!\geq\!5$, and 
(6) $Y_{\F_q((t))}(1,k)\!=\!\frac{q^k-1}{q-1}$ for any prime power $q$.  
Also: (7) $Y_{\Q_2}(1,k)\!\geq\!2k$, 
(8) $3\!\leq\!Y_{\Q_3}(1,2)\!\leq\!9$, (9) $Y_{\Q_p}(1,k)\!\geq\!2k-1$ 
for $p\!\geq\!3$, and (10) $Y_{\Q_p}(1,k)\!\leq\!k^2-k+1$ for 
$p\!>\!1+k$. \qed  
\end{thm} 
\begin{rem} 
\label{rem:low} The assertions above are immediate consequences  
of \cite[pg.\ 160]{descartes}, \cite[Cor.\ 2.1]{voorhoeve},
\cite[Example, pg.\ 286 \& pp.\ 289--290]{lenstra},
\cite[Thm.\ 1.4, Ex.\ 1.5, \& Thm.\ 1.6]{avendano2}, and
\cite[Sec.\ 2]{poonen}. Also, the polynomials 
$\prod^k_{i=1}(x_1-i)$, $3x^{10}_1+x^2_1-4$,\linebreak 
$x^{1+p^{p-1}}_1-(1+p^{p-1})x_1+p^{p-1}$, 
$\prod\limits_{z_1,\ldots,z_{k-1}\in\F_q}(x_1-z_1-z_2t-\cdots
-z_{k-1}t^{k-1})$, and $\prod^k_{i=1}(x^2_1-4^{i-1})$ respectively attain the 
number of roots stated in Assertions (1), (3), (5), (6), and (7). \dia 
\end{rem} 

\noindent 
$Y_L(1,1)$ can in fact grow without bound if we let $L$ range over 
arbitrary finite extensions of $\Q_p$.\footnote{For instance, when $L$ is the 
splitting field of $g(x_1)\!:=\!x^p_1-1$ over $\Q_p$, $g$ has roots 
$1,1+\mu_1,\ldots,1+\mu_{p-1}$ where the $\mu_i$
are distinct elements of $L$, each with valuation $\frac{1}{p-1}$
(see, e.g., \cite[pp.\ 102--109]{robert}).}  
Note also that for any local field $L\!\neq\!\C$ and fixed $(n,k)$, 
the supremum of the {\em total} number of roots of $F$ in $(L^*)^n$ --- with 
no restrictions on the phase of the coordinates --- 
is easily derivable from $Y_L(n,k)$ (see Proposition \ref{prop:iou} of Section
\ref{sec:wrap}). 

We treat the general multivariate case in Sections \ref{sub:simple} and \ref{sub:up}, where 
we state our main\linebreak 
\scalebox{.98}[1]{results. As a warm-up, let us first unite the simplest 
multivariate cases (proved in Section \ref{sec:wrap}).}   
\begin{prop}
\label{prop:simp} For any $k\!\leq\!0$, $n\!\geq\!1$, and any local field $L$,  
we have $Y_L(n,k)\!=\!0$. Also, $Y_{L}(n,1)\!=\!Y_L(1,1)^n$. In particular, 
$Y_{\Q_2}(n,1)\!=\!2^n$ and $Y_L(n,1)\!=\!1$ 
for all $L\!\in\!\{\C,\R\}\cup\{\Q_3,\Q_5,\ldots\} \cup\{\F_q((t))\; | \; q 
\text{ a prime power}\}$. 
\end{prop}  

\subsection{New, Simple Systems with Many Roots} 
\label{sub:simple} 
For any $j,N\!\in\!\N$ let $[j]_N\!\in\!\{0,\ldots,N-1\}$ denote the mod $N$ 
reduction of $j$. 
\begin{thm}  
\label{thm:big} 
For any local field $L$, 
$Y_L(n,2)\!\geq\!\max\left\{Y_L(1,1)^{n-1}Y_L(1,2),n+1\right\}$. 
More\linebreak 
\scalebox{.94}[1]{generally, $Y_L(n,k)\geq \max\left\{Y_L(1,1)^{n-k+1}
Y_L(1,2)^{k-1},Y_L\!\left(\left\lfloor \frac{n}{k-1}\right\rfloor,2 
\right)^{k-1-[n]_{k-1}}Y_L\!\left(\left\lfloor \frac{n}{k-1}\right
\rfloor+1,2\right)^{[n]_{k-1}}\right\}$}\linebreak 
when $n\!\geq\!k-1\!\geq\!1$, and  
$Y_L(n,k)\geq Y_L\!\left(1,\left\lfloor \frac{n+k-1}{n}\right\rfloor
\right)^{n-[k-1]_{n}}Y_L\!\left(1,\left\lfloor \frac{n+k-1}{n}
\right\rfloor+1\right)^{[k-1]_{n}}$ when $1\!\leq\!n\!\leq\!k-1$. More 
explicitly, the following lower bounds hold: 

\noindent 
\renewcommand{\arraystretch}{1.5}
\mbox{}\hfill
\begin{tabular}{c|c|c}
$L$ & $n\!\geq\!k-1\!\geq\!1$ & $1\!\leq\!n\!\leq\!k-1$ \\ \hline
$\R$ & $\left\lfloor \frac{n+k-1}{k-1}\right\rfloor
^{k-1-[n]_{k-1}}\left\lfloor \frac{n+2k-2}{k-1}\right
\rfloor^{[n]_{k-1}}$ & $\left\lfloor \frac{n+k-1}{n}\right\rfloor
^{n-[k-1]_n}\left\lfloor \frac{2n+k-1}{n}\right
\rfloor^{[k-1]_n}$ \\ 
$\Q_2$ & $2^n 3^{k-1}$ & $2^n \left\lfloor \frac{n+k-1}{n}\right\rfloor
^{n-[k-1]_n}\left\lfloor \frac{2n+k-1}{n}\right\rfloor^{[k-1]_n}$ \\ 
$\Q_p$ ($p\!\geq\!3$) & $\left\lfloor \frac{n+k-1}{k-1}\right\rfloor
^{k-1-[n]_{k-1}}\left\lfloor \frac{n+2k-2}{k-1}\right
\rfloor^{[n]_{k-1}}$ & \scalebox{.75}[1]
{$\left(2\left\lfloor \frac{n+k-1}{n}\right\rfloor-1
\right)^{n-[k-1]_n}\left(2\left\lfloor \frac{n+k-1}{n}\right\rfloor
+1\right)^{[k-1]_n}$} \\ 
$\F_q((t))$ & \scalebox{.7}[1]
{$\max\left\{q+1,\left\lfloor \frac{n+k-1}{k-1}\right\rfloor
\right\}^{k-1-[n]_{k-1}}\max\left\{q+1,\left\lfloor \frac{n+2k-2}{k-1}\right
\rfloor\right\}^{[n]_{k-1}}$} & \scalebox{.8}[1]
{$\left(\frac{q^{\left\lfloor \frac{n+k-1}{n}
\right\rfloor}-1}{q-1}\right)^{n-[k-1]_n}
\left(\frac{q^{\left\lfloor \frac{2n+k-1}{n}
\right\rfloor}-1}{q-1}\right)^{[k-1]_n}$}  
\end{tabular}
\hfill\mbox{} 
\end{thm} 

\noindent 
The lower bound $Y_\R(n,2)\!\geq\!n+1$ was first proved through an  
ingenious application of Dessins d'Enfants \cite{bihan}. We attain our 
more general lower bound for $Y_L(n,2)$ via an explicit family of polynomial 
systems instead. Note also that the $L\!=\!\R$ case of our general lower bound  
slightly improves an earlier 
$\left\lfloor \frac{n+k-1}{\min\{n,k-1\}}\right\rfloor^{\min\{n,k-1\}}$ 
lower bound from \cite{brs}. Non-trivial lower bounds, for 
$n\!\geq\!k-1\!\geq\!2$, were unknown for the non-Archimedean case. 

Letting $\Rn_+$ denote the
positive orthant, $\bL$ the algebraic closure of $L$, and defining\linebreak  
$\ord \; x\!:=\!-\log|x|$ in the Archimedean case, our new family of 
extremal systems can be described as follows:
\begin{thm} 
\label{thm:family} 
For any $n\!\geq\!2$, any local field $L$, and any $\eps\!\in\!L^*$ with 
generalized phase $1$ and $\ord\; \eps$ sufficiently large, the roots 
in $\bL^n$ of the $(n+2)$-nomial $n\times n$ system $G_\eps$ defined by 

\medskip 
\noindent 
\scalebox{.88}[1]{$\displaystyle{\left(x_1x_2-\left(\eps+x^2_1\right),
x_2x_3-\left(1+\eps x^2_1\right),
x_3x_4-\left(1+\eps^3x^2_1\right),
\ldots,x_{n-1}x_n-\left(1+\eps^{2n-5}x^2_1\right),
x_n-\left(1+\eps^{2n-3}x^2_1\right)\right)}$} 

\medskip 
\noindent
are all non-degenerate, lie in $(L^*)^n$, and have generalized phase $1$ for 
all their coordinates. In particular, $G_\eps$ has exactly 
$n+1$ non-degenerate roots in $\Rn_+$, $(\Q^*_p)^n$, or $(\F_q((t))^*)^n$ 
(each with generalized phase $1$ for all its coordinates), according as $\eps$ 
is $1/4$, $p$, or $t$.  
\end{thm} 

\noindent 
Explicit examples evincing $Y_\R(n,2)\!\geq\!n+1$ were previously known 
only for $n\!\leq\!3$ \cite{brs}. Our new extremal examples from 
Theorem \ref{thm:family} provide a new and arguably 
simpler proof that $Y_\R(n,2)\!\geq\!n+1$. 
We prove Theorems \ref{thm:big} and \ref{thm:family} in Sections  
\ref{sub:big} and \ref{sub:family}, respectively. 

\begin{rem} 
By construction, when we are over $\Q_p$ or $\F_q((t))$, the 
underlying tropical varieties of the zero sets defined by 
$G_\eps$ have a common form: they are each the Minkowski sum of 
an $(n-2)$-plane and a ``Y'' lying in a complementary $2$-plane. 
(See Section \ref{sec:mixed} for further background and Section 
\ref{sub:tropical} for some illustrations.) Furthermore, all these tropical 
varieties contain half-planes parallel to a single $(n-1)$-plane.  It is an 
amusing exercise to build such a collection of tropical varieties so that they 
have  at least $n+1$ isolated intersections. However, it is much more difficult 
to build a collection of polynomials whose 
tropical varieties have this property, and this constitutes 
a key subtlety behind Theorem \ref{thm:family}. \dia  
\end{rem} 

Another important construction underlying Theorem \ref{thm:family} is 
a particular structured family of univariate polynomials. 
\begin{lemma}  
\label{lemma:kait} For any $n\!\geq\!2$, the degree $n+1$ polynomial 
$R_n$ defined by 

\smallskip 
\noindent 
$u (1+\eps u)^2 (1+\eps^5 u)^2\cdots 
(1+\eps^{4\lfloor n/2\rfloor-3}u)^2 
- (\eps+u)^2 (1+\eps^3 u)^2 (1+\eps^7 u)^2 \cdots 
 (1+\eps^{4\lceil n/2\rceil-5}u)^2$ 

\smallskip 
\noindent 
has exactly $n+1$ roots in $\R_+$, $\Q^*_p$, or $\F_p((t))^*$, according as 
$\eps$ is $1/4$, $p$, or $t$. In particular, for these choices of $\eps$, 
all the roots of $R_n$ have generalized phase $1$.   
\end{lemma} 
\noindent 
We will see in Section \ref{sec:app} how the  
$R_n$ are part of a more general class of polynomials 
providing a bridge between fewnomial theory and algorithmic complexity. 
Lemma \ref{lemma:kait} is proved in Section \ref{sub:kait}.

\subsection{Upper Bounds: Known and Conjectural}  
\label{sub:up}  
That $Y_\R(n,k)\!<\!\infty$ for $n\!\geq\!2$ was first proved
around 1979 by Khovanskii and Sevastyanov \cite{kho,few}, yielding
an explicit, singly-exponential upper bound.
Based on the seminal results \cite[Pg.\ 105]{vandenef} and
\cite[Thm.\ 2]{lipshitz}
the second author proved in \cite[Thm.\ 1]{fin} that
$Y_{L}(n,k)\!<\!\infty$ for any fixed
$n$, $k$, and non-Archimedean field $L$ of characteristic zero.
(See \cite{amd} and the table below for explicit upper bounds.)
The finiteness of $Y_{\F_q((t))}(n,k)$ for $n\!\geq\!2$ remains unknown, in
spite of recent results of Avenda\~no and Ibrahim \cite{ai2} giving explicit
upper bounds for the number of roots in $L^n$ of a large class of $n\times n$
systems over any non-Archimedean local field $L$.
 
We will use Landau's $O$-notation for asymptotic upper bounds 
modulo a constant multiple, along with the companion $\Omega$-notation 
for asymptotic lower bounds.  
The best known upper and lower bounds on 
$Y_L(n,k)$ (as of November 2012), for $L\!\in\!\{\R,\Q_3,\Q_5,\ldots\}$ 
and $n,k\!\geq\!2$, are the following: 

\renewcommand{\arraystretch}{2}
\medskip 
\noindent 
\mbox{}\hfill
\scalebox{1}[1]{
\begin{tabular}{c|c|c}
$L$ & Upper Bound on $Y_L(n,k)$ & Lower Bound on $Y_L(n,k)$ \\ \hline 
$\R$ & 
$2^{O(k^2)}n^{k-1}$ \hspace{.7cm} \cite{bs}$^4$ & 
$\Omega\!\left(\left\lfloor\frac{n+k-1}{\min\{n,k-1\}}\right\rfloor
\right)^{\min\{n,k-1\}}$  
(Theorem \ref{thm:big} here) \\ 
$\Q_p$ & $\left(O\!\left(k^3 n \log k\right)\right)^n$ 
\cite{amd} & $\Omega\!\left(\left\lfloor\frac{n+k-1}{\min\{n,k-1\}}\right 
\rfloor\right)^{\min\{n,k-1\}}$ (Theorem \ref{thm:big} here) 
\end{tabular}}\hfill\mbox{} 
\addtocounter{footnote}{1}
\footnotetext{
While there have been important recent refinements to this bound (e.g., 
\cite{rss}) the asymptotics of \cite{bs} have not yet been 
improved in complete generality.} 

\renewcommand{\arraystretch}{2}

\medskip 
\noindent 
Also, Bertrand, Bihan, and Sottile proved 
the (tight) upper bound $Y_\R(n,2)\!\leq\!n+1$ in  
\cite{bbs}. The implied $\Omega$-constants 
above can be taken to be $1$. 

Most importantly, note that for the Archimedean case (resp.\ the $p$-adic 
rational case with $p\!\geq\!3$), $Y_L(n,k)$ is bounded from above by a  
polynomial in $n$ when 
$k$ is fixed (resp.\ a polynomial in $k$ when $n$ is fixed). Based on this 
asymmetry of upper bounds, the second author posed the following conjecture 
(mildly paraphrased) at his March 20 Geometry Seminar talk at the Courant 
Institute in March 2007. 
\begin{localconj} \mbox{}\\ 
There are absolute constants $C_2\!\geq\!C_1\!>\!0$ such that, 
for any $L\!\in\!\{\C,\R,\Q_3,\Q_5,\ldots\}$ 
and any $n,k\!\geq\!2$, we 
have $(n+k-1)^{C_1\min\{n,k-1\}}\leq Y_L(n,k)\leq 
(n+k-1)^{C_2\min\{n,k-1\}}$. 
\end{localconj} 
\begin{rem}
Should the Local Fewnomial Conjecture be true, it is likely that
similar bounds can be asserted for the number of roots counting multiplicity, 
in the characteristic zero case.
This is already known for $(L,n)\!=\!(\R,1)$ \cite{wang}, and
\cite{lenstra,amd} provide evidence for the $p$-adic rational case.
Note, however, that the equality $(x_1+1)^{q^m+1}
\!=\!x^{q^m+1}_1+x^{q^m}_1+x_1+1$ over $\F_q$ (as observed in \cite{poonen}) 
tells us that for $L$ of positive characteristic it is impossible to 
count roots over $L^*$ --- {\em with multiplicity} --- solely as a function 
of $n$, $k$, and $L$. \dia
\end{rem}

\noindent 
Theorem \ref{thm:big} thus reveals the lower bound 
of the Local Fewnomial Conjecture to be true (with $C_1\!=\!1$) 
for the special case $k\!=\!2$.  
From our table above we also see that the upper bound 
from the Local Fewnomial Conjecture holds for $n\!\leq\!k-1$ 
(at least for $C_2\!\geq\!7$), in the $p$-adic rational  
setting. We intend for 
our techniques here to be a first step toward establishing the Local Fewnomial  
Conjecture for $n\!>\!k-1$ in the $p$-adic rational setting. 

Note that the maximal number of roots in $\Csn$ 
of an $(n+k)$-nomial $n\times n$ system $F$ over $\C$ is undefined 
for any fixed $n$ and $k$: consider $((x^d_1-1)\cdots(x^d_1-k),x_2-1,\ldots,
x_n-1)$ as $d\longrightarrow\infty$. Nevertheless, the maximal number of roots 
in $\Rn_+$ is well-defined and finite for any fixed $n,k\!\geq\!1$. The latter 
assertion is a very special case of Khovanski's {\em Theorem on Complex 
Fewnomials} (see \cite[Thm.\ 1 (pp.\ 82--83), Thm.\ 2 (pp.\ 87--88), and 
Cor.\ 3$'$ (pg.\ 88)]{few}), which estimates the number of roots in angular 
sub-regions of $\Cn$ for a broad class of analytic functions. 
\cite{few} does not appear to state any explicit upper bounds for 
$Y_\C(n,k)$, but one can in fact show (see Section 
\ref{sec:wrap}) that it suffices to study the real case. 
\begin{thm} 
\label{thm:hah} 
For all $n,k\!\geq\!1$, we have $Y_\C(n,k)\!=\!Y_\R(n,k)$. 
\end{thm}

We now discuss the number of roots, over a local field, of certain non-sparse 
univariate polynomials that nevertheless admit a 
compact expression, e.g., $(x^9_1+1)^{1000}-(x_1-3)^{2^8}$. This refinement 
leads us to computational number theory and variants of the 
famous {\em $\pp$ vs.\ $\np$ Problem}. As we will see shortly, complexity 
theory leads us to challenging open problems that can be stated entirely 
within the context of arithmetic geometry. 

\section{Applications and New Conjectures on Straight-Line Programs} 
\label{sec:app} 
To better discuss the connections between structured polynomials and 
algorithms let us first introduce the notions of {\em input size} and {\em 
complexity} through a concrete example. \cite{bas} is an outstanding 
reference for basic algorithmic number theory and \cite{sipser,papa,cxity,
fortnow,liptonblog} are among many excellent sources for further background on 
complexity theory and the history of the $\pp$ vs.\ $\np$ Problem.  
\begin{ex} 
Consider the following problem: 
\begin{quote} 
{\bf A:} Given any prime $p$ and $f\!\in\!\F_p[x_1]$ with degree $d$ 
and $d\!<\!p$, decide whether\linebreak 
\mbox{}\hspace{.65cm}$f$ has a root in $\F_p$. 
\end{quote} 
Let us naturally define the input size of an 
instance $(p,f)$ of Problem A as the number of decimal digits needed to write 
down $p$ and the monomial term expansion of $f$. (Thus, for example, 
$a+bx^{11}+cx^d$ would have size $O(\log p)$ since $a,b,c,d\!\in\!\{0,
\ldots,p-1\}$.) To measure the complexity of 
a computation over $\F_p$, we can then simply count the number of digit by 
digit operations (i.e., addition, subtraction, multiplication, and parity 
checking) that we use. For instance, via fast mod $n$ arithmetic (e,g., 
\cite[Ch.\ 5]{bas}), it is easy to see that evaluating $f$ at a point in 
$\F_p$ has complexity near-linear in the input size (a.k.a.\ near-linear 
{\em time}). 

Curiously, no method with complexity polynomial in the input size is 
known for Problem A, although a putative root can be certainly be 
verified in polynomial-time.\footnote{Technically, Problem A is in 
$\np$, and is $\np$-hard with respect to randomized reductions 
\cite{bcr}.} \dia  
\end{ex} 

The complexity of evaluating a polynomial turns out to be a more intrinsic 
measure of its size than counting digits in monomial term expansions. In 
particular, many non-sparse polynomials can still be evaluated efficiently 
since they may admit other kinds of compact expressions. 
One central notion refining our preceding definition of input size 
is {\it straight-line program (SLP) complexity}.
\begin{dfn} 
For any field $K$ and $f\!\in\!K[x_1]$ let $s(f)$ --- the {\em SLP 
complexity of $f$} --- denote the smallest $n$ 
such that $f\!=\!f_n$ identically where the sequence $(f_{-N},\ldots,f_{-1},
f_0,\ldots,f_n)$ satisfies the following conditions: 
$f_{-1},\ldots,f_{-N}\!\in\!K$, $f_0\!:=\!x_1$, and, for all $i\!\geq\!1$,  
$f_i$ is a sum, difference, or product of some pair of elements 
$(f_j,f_k)$ with $j,k\!<\!i$. Finally, for any $f\!\in\!\Z[x_1]$, we 
let $\tau(f)$ denote the obvious analogue of $s(f)$ where the definition 
is further restricted by assuming $N\!=\!1$ and $f_{-1}\!:=\!1$. \dia 
\end{dfn} 

\noindent 
Note that we always have $s(f)\!\leq\!\tau(f)$ since $s$ does not 
count the cost of computing large integers (or any constants). 
\begin{ex} 
Evaluating $x^{2^k}_1$ via recursive squaring (i.e., 
$(\cdots(x^2_1)^2\cdots)^2$), and employing the binary 
expansion of $d$, it is easily checked that 
$s\!\left(x^d_1\right)\!=\!\tau\!\left(x^d_1\right)\!=\!O\!\left(\log^2 
d\right)$. One in fact has $\tau(n)\!\leq\!2\log_2 n$ for any 
$n\!\in\!\N$ \cite[Prop.\ 1]{svaiter} and, when $n$ is a difference of two 
nonnegative integers with at most $\delta$ nonzero digits in their binary 
expansions, we also obtain $s(n)\!=\!1$ and 
$\tau(n)\!=\!O\!\left(\delta(\log\log|n|)^2\right)$. 
See also \cite{brauer,moreira} for further background. \dia   
\end{ex} 

Relating SLP complexity to the number of rational 
roots of polynomials provides a delightfully direct way to go from the theory 
of sparse polynomials to deep open questions in complexity theory and 
computational number theory. In what follows, we let $Z_R(f)$ denote the set 
of roots of $f$ in a ring $R$, and use $\#S$ for the cardinality of a set $S$.  
\begin{thm} 
\label{thm:tau} \mbox{}\\  
I. (See \cite[Thm.\ 3, Pg.\ 127]{bcss} and \cite[Thm.\ 1.1]{burgtau}.)
Suppose that for all nonzero $f$\linebreak 
\mbox{}\hspace{.5cm}$\in\!\Z[x_1]$ we have 
$\#Z_\Z(f)\!\leq\!(\tau(f)+1)^{O(1)}$.
Then $\pp_\C\!\neq\!\np_\C$, and the permanent of $n\times n$ 
ma-\linebreak 
\mbox{}\hspace{.5cm}trices cannot be computed by constant-free, 
division-free arithmetic circuits of size $n^{O(1)}$. 

\noindent 
II. (Weak inverse to (I) \cite{lipton}.\footnote{
Lipton's main result from \cite{lipton} is in fact stronger, allowing for 
rational roots and primes with a mildly differing number of digits.}) If 
there is an $\eps\!>\!0$ and a sequence $(f_n)_{n\in\N}$ 
of polynomials\linebreak 
\mbox{}\hspace{.7cm}in $\Z[x_1]$ satisfying:\\ 
\mbox{}\hspace{1.2cm}(a) $\#Z_\Z(f_n)\!>\!e^{\tau(f_n)^\eps}$ for all 
$n\!\geq\!1$ \ and  \
(b) $\deg f_n, \max\limits_{\zeta\in Z_\Z(f)}|\zeta|\!\leq\!
2^{(\log \#Z_\Z(f_n))^{O(1)}}$\\  
\mbox{}\hspace{.7cm}then, for infinitely many $n$, at least 
$\frac{1}{n^{O(1)}}$ of the $n$ digit integers that are products of 
exactly\linebreak 
\mbox{}\hspace{.7cm}two distinct primes (with an equal number of 
digits) can be factored by a Boolean circuit\linebreak
\mbox{}\hspace{.7cm}of size $n^{O(1)}$.  

\noindent
III. (Number field analogue of (I) implies Uniform Boundedness \cite{cheng}.) 
Suppose that for\linebreak  
\mbox{}\hspace{.8cm}any number field $K$ and $f\!\in\!K[x_1]$ we have  
$\#Z_K(f)\!\leq\!c_1 1.0096^{s(f)}$, with $c_1$ depending only\linebreak   
\mbox{}\hspace{.8cm}on $[K:\Q]$. Then there is a constant $c_2\!\in\!\N$ 
depending only on $[K:\Q]$ such that for any\linebreak 
\mbox{}\hspace{.8cm}elliptic curve $E$ over $K$, the torsion 
subgroup of $E(K)$ has order at most $c_2$. \qed 
\end{thm} 

\noindent 
The hypothesis in Part (I) is known as the {\em (Shub-Smale) 
$\tau$-Conjecture}, and was also stated as the fourth problem on Smale's 
list of the most important problems for the $21\stst$ century \cite{21a,21b}. 
Mike Shub informed the authors in late 2011 that,  
should the $\tau$-Conjecture hold, its $O$-constant should be at least 
$2$. The complexity classes $\pp_\C$ and $\np_\C$ are respective analogues 
(for the BSS model over $\C$ \cite{bcss}) of the well-known complexity 
classes $\pp$ and $\np$. 
(Just as in the famous $\pp$ vs.\ $\np$ Problem, the equality of $\pp_\C$ and 
$\np_\C$ remains an open question.) 
The assertion on the hardness of the permanent in Theorem \ref{thm:tau}
is also an open problem and its proof would be a major step toward
solving the {\em $\vp$ vs.\ $\vnp$ Problem} ---
Valiant's algebraic circuit analogue of the $\pp$ vs.\ $\np$ Problem
\cite{valiant,burgcook,koiran,jml}.

The hypothesis of Part (II) merely posits a sequence of polynomials 
violating the $\tau$-Conjecture in a weakly exponential manner. 
The conclusion in Part (II) would violate a widely-believed 
version of the cryptographic hardness of integer factorization. 

Some evidence toward the hypothesis of Part (III) is provided by 
\cite[Thm.\ 1]{add}, which gives the upper bound 
$\#Z_K(f)\!\leq\!2^{O(\sigma(f)\log\sigma(f))}$. The quantity $\sigma(f)$ is 
the {\em additive} 
complexity of $f$ \cite{grigoadd,add} and is bounded from above by $s(f)$.   
The conclusion in Part (III) is the famous {\em Uniform Boundedness 
Theorem}, due to Merel \cite{merel}. 
Cheng's conditional proof (see 
\cite[Sec.\ 5]{cheng}) is dramatically simpler and would yield effective bounds 
significantly improving known results (e.g., those of Parent \cite{parent}). 
In particular, the $K\!=\!\Q$ case of the hypothesis of Part (III) 
would yield a new proof (less than a page 
long) of Mazur's landmark result on torsion points \cite{mazur}. 

A natural approach to the $\tau$-Conjecture would be to 
broaden it to inspire a new set of techniques, or rule out overly 
optimistic extensions. For instance, 
one might suspect that the number of roots of $f$ in a field $L$ containing 
$\Z$ could also be polynomial in $\tau(f)$, thus allowing us to consider  
techniques applicable to $L$. 
For $L$ a number field, the truth of such an  
extension of the $\tau$-Conjecture expands its implications into 
arithmetic geometry, as we already saw in Part (III) of Theorem \ref{thm:tau}. 
However, the truth of any global field analogue of the 
$\tau$-Conjecture remains unknown. 

Over local fields, we now know that the most naive extensions 
break down quickly: There are well-known examples $(f_n)_{n\in\N}$,  
from the dynamical systems and algorithms literature, 
with $\tau(f_n)\!=\!O(n)$ and $f_n$ having $2^n$ real roots (see, e.g., 
\cite{bocook,perrucci}). 
Constructing such ``small but mighty'' polynomials over $\Q_p$ is also 
possible, even over several such fields at once. 
\begin{ex} 
\label{ex:slp} 
Let $S$ be any non-empty finite set of primes, $c_S\!:=\!\prod_{p\in S}p$,  
$k\!:=\!\max S$, and consider the recurrence satisfying 
$h_1\!:=\!x_1(1-x_1)$ 
and $h_{n+1}\!:=\!\left(c^{3^{n-1}}_S-h_{n}\right)h_n$ for all 
$n\!\geq\!1$. Then $\frac{h_{n}(x_1)}{x_1(1-x_1)}\!\in\!\Z[x_1]$ has 
degree $2^n-2$, exactly $2^n-2$ roots in $\Z_p$ for each 
$p\!\in\!S$, and $\tau\!\left(\frac{h_{n}(x_1)}{x_1(1-x_1)}
\right)=O(n+\#S\log k)$. However, $\frac{h_{n}(x_1)}
{x_1(1-x_1)}$ has no real roots, and thus no integer roots. 
(Proofs of these facts are provided in Section \ref{sub:slp}.) \dia  
\end{ex}

To the best of our knowledge, the $\tau$-Conjecture still has no 
counter-examples. Indeed, all known families of ``small but mighty''  
polynomials are of a very particular recursive form, and have few 
(if any) integer roots at all.  
So let us now formulate a potentially safer extension of the $\tau$-Conjecture 
to local fields, and apply it to a more restricted family of expressions: 
{\em sum-product-sum (SPS) polynomials}. 
\begin{dfn} 
(See \cite[Sec.\ 3]{koiran}.) 
Let us define $\sps(k,m,t,d,\delta)$ to be the family of non-\linebreak  
constant polynomials presented in the form 
$\sum^k_{i=1} \prod^m_{j=1} f_{i,j}$ where, for all $i$ and $j$, \\
\mbox{}\hspace{1cm}(1) $f_{i,j}\!\in\!\Z[x_1]\!\setminus\!\{0\}$ has degree 
$\leq\!d$ and $\leq\!t$ monomial terms\\  
\mbox{}\hspace{1cm}(2) each coefficient of $f_{i,j}$ has absolute value 
$\leq\!2^d$, and is the difference of two nonneg-\\ 
\mbox{}\hspace{1.65cm}ative integers with at 
most $\delta$ nonzero digits in their binary expansions. \dia 
\end{dfn}

\noindent 
For instance, it is easily checked that the univariate polynomial\\  
\mbox{}\hfill $\left(7y^{97139}_1-9y^7\right)\left(24y^{45}_1+1000y^{131}_1
\right)+y^{99}_1$\hfill\mbox{}\\ 
lies in $\sps(2,2,2,97139,2)$. The family $\sps(k,m,t,d,\delta)$ is motivated 
by recent advances in circuit complexity \cite{agvinay,koiran}. 
SPS polynomials have also (implicitly) appeared earlier in fewnomial theory: 
\cite[Lemma 2]{tri}, \cite[Prop.\ 4.2, pg.\ 375]{bbs}, and \cite[Thm.\ 1]
{avendano}, in rather different notation,  
respectively derived upper bounds on the number of real roots of certain 
sub-families of $\sps(k,m,2,1,\delta)$, $\sps(2,m,d+1,d,\delta)$, 
and $\sps(k,2,2,1,\delta)$, independent of $\delta$.
Noting that $\tau(f)\!=\!(kmt+\delta+\log d)^{O(1)}$ for any 
$f\!\in\!\sps(k,m,t,d,\delta)$, we see that the following recent result 
of Koiran significantly strengthens part of Assertion (I) of Theorem 
\ref{thm:tau}. 
\begin{thm}  
\label{thm:koi} 
\cite[Conj.\ 1]{koiran}  
Suppose that for all $k,m,t,d,\delta\!\in\!\N$ and 
$f\!\in\!\sps(k,m,t,d,\delta)$, we have 
$\#Z_\Z(f)\!=\!(kmt+\delta+\log d)^{O(1)}$.  
Then the permanent of $n\times n$ 
matrices cannot be computed by 
constant-free, division-free arithmetic circuits of size $n^{O(1)}$. \qed  
\end{thm} 

In \cite{koiran}, Koiran suggests further that the number of real 
roots may also satisfy a bound like the one above. We propose a more 
flexible conjecture. 
\begin{adelic} 
For any $k,m,t,d,\delta\!\in\!\N$ and $f\!\in\!\sps(k,m,t,d,\delta)$, 
there is a field\linebreak  
\scalebox{.97}[1]{$L\!\in\!\{\R,\Q_2,\Q_3,\Q_5,\dots\}$ such that 
$f$ has no more than $(kmt+\delta+\log d)^{O(1)}$ distinct roots in $L$.}  
\end{adelic} 

\noindent 
The Adelic $\tau$-Conjecture clearly implies the 
hypothesis of Theorem \ref{thm:koi}. (Some evidence 
toward the Adelic $\tau$-Conjecture appears in \cite{unreal}.)  
So we pose our conjecture mainly to advocate 
adding $p$-adic techniques to the real-analytic toolbox put forth in 
\cite[Sec.\ 6]{koiran} and \cite{koiwron}.

\section{Background: From Triangles to Toric Deformations and Tropical 
Varieties}
\label{sec:mixed} 
Our first step toward building systems with maximally many roots is a 
polyhedral\linebreak 
construction (Lemma \ref{lemma:tri} below) with several useful 
algebraic consequences. We refer the reader to the excellent 
book \cite{triang} for further background on triangulations and liftings. 

Let $\conv \cA$ denote the convex hull of any set $\cA\!\subseteq\!\Rn$.  
Assuming $\cA$ is finite, we say that a triangulation of $\cA$ is 
{\em coherent} (or {\em regular}) 
iff its simplices are exactly the domains of linearity for some function 
$\ell : \conv\cA \longrightarrow \R$ 
that is convex, continuous, and piecewise linear. (For $n\!\geq\!2$ 
and $\#\cA\!\geq\!6$ one can easily find non-coherent triangulations 
\cite{triang}.) We call $\ell$ 
a {\em lifting} of $\cA$ (or a lifting of $\conv \cA$), and 
we let $\hA\!:=\!\{(a,\ell(a))\; | \; a\!\in\!\cA\}$. Abusing 
notation slightly, we also refer to $\hA$ as a {\em lifting of $\cA$ 
(with respect to $\ell$)}. 
\begin{rem} 
It follows directly from our last definition that a lifting 
function $\ell$ on $\conv \cA$ is uniquely determined by the values of 
$\ell$ on $\cA$. So we will henceforth specify such $\ell$  
by specifying just the restricted image $\ell(\cA)$. \dia 
\end{rem} 
Recall also that $\supp(f)$ denotes the set of exponent vectors
(a.k.a.\ the {\em support} or\linebreak {\em spectrum}) of $f$.
\begin{ex}
Consider $f(x)\!:=\!1-x_1-x_2+\frac{6}{5}(x^4_1x_2+x_1x^4_2)$. Then
$\supp(f)\!=\!\{(0,0),(1,0),$\linebreak
$(0,1),(1,4),(4,1)\}$ and has convex hull a pentagon.
It is then easily checked that there are exactly $5$ possible 
triangulations for $\supp(f)$, all of which happen to be coherent:\\
\mbox{}\hfill\epsfig{file=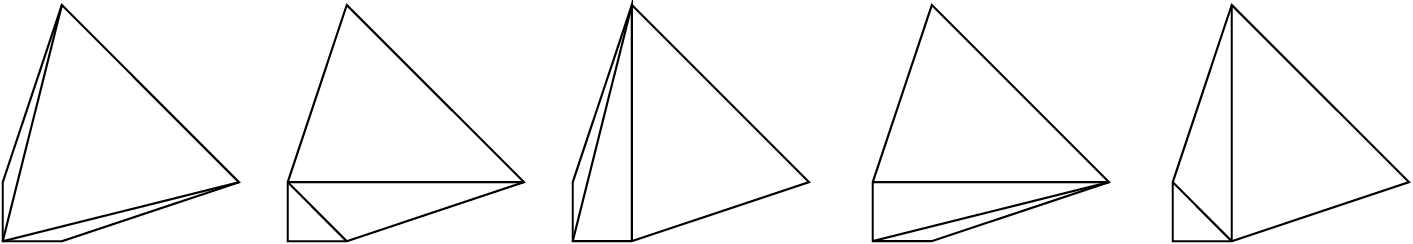,width=.7\textwidth} \hfill \dia
\end{ex}
\begin{dfn} 
\label{dfn:lower}  
(See also \cite{hs}.) 
For any polytope 
$\hat{Q}\!\subset\!\R^{n+1}$, we call a face $\hat{P}$ of $\hat{Q}$ a 
{\em lower} face iff $\hat{P}$ has an inner normal with positive 
$(n+1)\stst$ coordinate. 
Letting $\pi : \R^{n+1} \longrightarrow \Rn$ denote the natural projection 
forgetting the last coordinate, the lower facets of $\hQ$ thus induce a 
natural polyhedral subdivision $\Sigma$ of $Q\!:=\!\pi\!\left(\hat{Q}\right)$. 
In particular, if $\hat{Q}\!\subset\!\R^{n+1}$ is a Minkowski sum of the form 
$\hat{Q}_1+\cdots+\hat{Q}_n$ where the 
$\hat{Q}_i$ are polytopes of dimension $\leq n+1$, $\hat{E}_i$ is a lower edge 
of $\hQ_i$ for all $i$, and $\hP\!=\!\hE_1+\cdots+\hE_n$ is a 
lower facet of $\hQ$, then we call $\hP$ a {\em mixed} lower facet 
of $\hQ$. 
Also, the resulting cell $\pi\!\left(\hP\right)\!=\!\pi\!\left(\hE_1\right)+
\cdots+\pi\!\left(\hE_n\right)$ of $\Sigma$ is called a 
{\em mixed cell} of $\Sigma$. \dia  
\end{dfn} 
\begin{ex}
\label{ex:g2} 
Let us consider the family of systems $G_\eps$ from Theorem \ref{thm:family} 
for $n\!=\!2$. In particular, let $(\cA_1,\cA_2)$ be the pair of supports of 
$G_\eps$, and let $(Q_1,Q_2)$ be the corresponding pair of convex hulls in 
$\R^2$. Let us also define a pair of liftings $(\ell_1,\ell_2)$ via the
exponents of the powers of $\eps$ appearing in the corresponding monomial
terms. More precisely, $\ell_1$ sends $(0,0)$, $(2,0)$, and $(1,1)$ 
respectively to $1$, $0$, and $1$; and $\ell_2$ sends 
$(1,1)$, $(2,0)$, and $(0,1)$ respectively to $0$, $1$, and $0$. 
These lifting functions then affect the shape of the lower hull of the 
Minkowski sum $\hQ_1+\hQ_2$ of lifted polygons, which in turn fixes
a subdivision $\Sigma_{\ell_1,\ell_2}$ of $Q_1+Q_2$ via the images of the 
lower facets of $\hQ_1+\hQ_2$ under $\pi$. (See the illustration below.)
The mixed cells of $\Sigma_{\ell_1,\ell_2}$,\linebreak  

\vspace{-.8cm}
\noindent
\begin{minipage}[t]{.4\textwidth}
\vspace{0pt}
\mbox{}\hspace{1cm}\epsfig{file=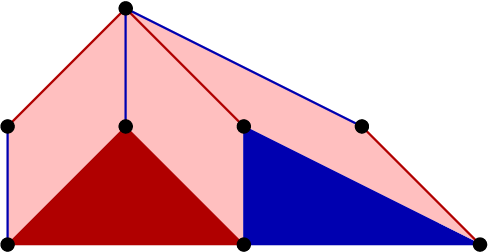,height=.9in}
\end{minipage}
\begin{minipage}[t]{.6\textwidth}
\vspace{0pt}
for this particular lifting, correspond to the lighter (pink) parallelograms: 
from left to right, they are exactly $E_{1,0}+E_{2,0}$, 
$E_{1,1}+E_{2,0}$, and $E_{1,1}+E_{2,1}$, where $E_{1,s}$ (resp.\ $E_{2,s}$) 
is an edge of $Q_1$ (resp.\ $Q_2$) for all 
$s$. More precisely, $E_{1,0}$, $E_{1,1}$, $E_{2,0}$,\linebreak  
\end{minipage} 

\vspace{-.3cm} 
\noindent 
and $E_{2,1}$ are respectively the convex hulls of $\{(0,0),(1,1)\}$, 
$\{(1,1),(2,0)\}$, 
$\{(0,0),(0,1)\}$, and $\{(0,1),(2,0)\}$. Note also that these mixed cells,   
through their expression as edges sums (and the obvious correspondence between 
vertices and monomial terms), correspond naturally to three binomial systems. 
In order, they are $(x_1x_2-\eps,x_2-1)$ , $(x_1x_2-x^2_1,x_2-1)$, and 
$(x_1x_2-x^2_1,x_2-\eps x^2_1)$. In particular, the first (resp.\ second) 
polynomial of each such pair is a sub-sum of the first (resp.\ second) 
polynomial of $G_\eps$. \dia 
\end{ex}
\begin{dfn} 
\label{dfn:mixed}  
(See also \cite{hs,ewald,why}.) 
Let $\cA_1,\ldots,\cA_n\!\subset\!\Rn$ be finite point sets 
with respective convex hulls $Q_1,\ldots,Q_n$. 
Also let  $\ell_1,\ldots,\ell_n$ be respective lifting functions for 
$\cA_1,\ldots,\cA_n$ and consider the  
polyhedral subdivision $\Sigma_{\ell_1,\ldots,\ell_n}$ of 
$Q\!:=\!Q_1+\cdots+Q_n$ obtained via the images of the lower facets
of $\hQ$ under $\pi$. In particular, if 
$\dim \hP_1+\cdots+\dim \hP_n\!=\!n$ for every lower facet of $\hQ$ 
of the form $\hP_1+\cdots+\hP_n$, then we say that $(\ell_1,
\ldots,\ell_n)$ is {\em mixed}. For any mixed $n$-tuple of 
liftings we then define the {\em mixed volume of 
$(Q_1,\ldots,Q_n)$} to be $\cM(Q_1,\ldots,Q_n):=\!\!\!\!\!\!\!\sum
\limits_{\substack{C \text{ a mixed cell}\\ \text{of } 
\Sigma_{\ell_1,\ldots,\ell_n}}} \!\!\!\!\! \vol(C)$,  
following the notation of Definition \ref{dfn:lower}. \dia 
\end{dfn} 

\noindent 
As an example, the mixed volume of the two triangles from Example 
\ref{ex:g2}, relative to the stated (mixed) lifting, is the sum of the areas 
of the three parallelograms in the illustration, i.e., $3$. 
\begin{thm} \label{thm:mixed} 
(See \cite[Ch.\ IV, pg.\ 126]{ewald} and \cite{hs}.)  
The formula for $\cM(Q_1,\ldots,Q_n)$ from Definition \ref{dfn:mixed} 
is independent of the underlying mixed $n$-tuple of\linebreak liftings 
$(\ell_1, \ldots,\ell_n)$. Furthermore, if $Q'_1,\ldots,Q'_n\!\subseteq\!\Rn$ 
are any polytopes with $Q'_i\!\supseteq\!Q_i$ for all $i$, then 
$\cM(Q_1,\ldots,Q_n)\!\leq\!\cM(Q'_1,\ldots,Q'_n)$. Finally, the 
$n$-dimensional mixed volume\linebreak 
satisfies $\cM(Q,\ldots,Q)\!=\!n!\vol(Q)$ for any polytope $Q\!\subset\!\Rn$. 
\qed 
\end{thm} 
\begin{lemma} 
\label{lemma:tri} 
Let $n\!\geq\!2$, and let $\bO$ and $e_i$ respectively denote the origin and 
$i\thth$ standard basis vector in $\R^{n+1}$. Consider the triangles 
$\hT_1\!:=\!\conv\{e_{n+1},2e_1,e_1+e_2\}$,\linebreak  
$\hT_n\!:=\!\conv\{\bO,2e_1+(2n-3)e_{n+1},e_n\}$, and 
$\hT_i\!:=\conv\{\bO,2e_1+(2i-3)e_{n+1},e_i+e_{i+1}\}$ for all 
$i\!\in\!\{2,\ldots,n-1\}$. 
Then the Minkowski sum $\hT\!:=\!\hT_1+\cdots+\hT_n$ has exactly 
$n+1$ mixed lower facets. More precisely, 
for any $j\!\in\!\{0,\ldots,n\}$, we can obtain a unique mixed lower facet, 
$\hP_j:=\hE_{1,1}+\cdots+\hE_{j,1}+\hE_{j+1,0}+\cdots
+\hE_{n,0}$, with $\vol\!\left(\pi\!\left(\hP_j\right)\right)\!=\!1$, 
in the following manner: for all $i\!\in\!\{1,\ldots,n\}$, define $\hE_{i,1}$ 
(resp.\ $\hE_{i,0}$) to be the convex hull of the second (resp.\ first) and 
third listed vertices for $\hT_i$. Finally, $\cM\!\left(
\pi\!\left(\hT_1\right),\ldots,\pi\!\left(\hT_n\right)\right)\!=\!n+1$ and, 
for each $j\!\in\!\{0,\ldots,n\}$, 
the vector $v_j\!:=\!e_{n+1}+e_1-\sum^{j}_{i=1} (j+1-i)e_i$ 
is a nonzero inner normal for the lower facet $\hP_j$. 
\end{lemma} 

\noindent 
Lemma \ref{lemma:tri} is our key polyhedral result and is proved in Section 
\ref{sub:tri} and illustrated in Example \ref{ex:3d} below. 

The next result we need is a beautiful generalization, by Bernd Sturmfels, 
of {\em Viro's}\linebreak 
{\em Theorem}. We use $\partial Q$ for the boundary of a polytope $Q$.   
\begin{dfn}
\label{dfn:viro} 
Suppose $\cA\!\subset\!\Zn$ is finite and $\vol(\conv \cA)\!>\!0$. 
We call any function $s : \cA \longrightarrow \{\pm\}$ a
{\em distribution of signs for $\cA$}, and we call any pair  
$(\Sigma,s)$ with $\Sigma$ a coherent triangulation of $\cA$ 
a {\em signed} (coherent) triangulation of $\cA$. We also 
call any edge of $\Sigma$ with vertices of opposite sign an 
{\em alternating edge}. 

Given a signed triangulation for $\cA$ we then define a piece-wise linear 
manifold  --- the {\em Viro diagram} $\cV_\cA(\Sigma,s)$ --- in the 
following local manner: For any $n$-cell $C\!\in\!\Sigma$,
let $L_C$ be the convex hull of the set of midpoints of the alternating edges 
of $C$, and then define $\cV_\cA(\Sigma,s)\!:=\!\bigcup\limits_{\substack{C 
\text{ an } n\text{-cell}\\ \text{of } \Sigma}} L_C\setminus\partial 
\conv(\cA)$. Finally, when $\cA\!=\!\supp(f)$ and 
$s$ is the corresponding sequence of coefficient signs, then we call 
$\cV_{\Sigma}(f)\!:=\!\cV_\cA(\Sigma,s)$ the {\em Viro diagram of $f$}. \dia
\end{dfn}

Viro's Theorem (see, e.g., Proposition 5.2 and Theorem 5.6 of
\cite[Ch.\ 5, pp.\ 378--393]{gkz94} or \cite{viro}) states that, under 
certain conditions, one may 
find a triangulation $\Sigma$ with the positive zero set of $f$ 
homeomorphic to $\cV_\Sigma(f)$. {\em Sturmfels' Theorem for Complete 
Intersections} \cite[Thm.\ 4]{sturmfels} extends this to polynomial systems, 
and we will need just the $n\times n$ case. 
\begin{dfn} 
\label{dfn:alt} 
Suppose $\cA_1,\ldots,\cA_n\!\subset\!\Zn$ and each 
$\cA_i$ is endowed with a lifting $\ell_i$ and a 
distribution of signs $s_i$. Then, following the notation of Definition 
\ref{dfn:mixed}, we call a mixed cell $E_1+\cdots+E_n$ of 
$\Sigma_{\ell_1,\ldots,\ell_n}$ an {\em alternating mixed cell of 
$(\Sigma_{\ell_1,\ldots,\ell_n},s_1,\ldots,s_n)$} iff each edge $E_i$ is 
alternating (as an edge of the triangulation of $\cA_i$ induced by 
$\ell_i$). \dia 
\end{dfn} 
\begin{ex} 
Returning to Example \ref{ex:g2}, it is clear that, when $\eps\!\in\!\Rs$, 
we can endow the supports of $G_\eps$ with the distribution of signs 
corresponding to the underlying coefficients. In particular, when 
$\eps\!>\!0$, each of the $3$ mixed cells is alternating. \dia 
\end{ex} 

\begin{sturmf} 
Suppose $\cA_1,\ldots,\cA_n$ are finite subsets of $\Zn$, 
$(c_{i,a}\; | \; i\!\in\!\{1,\ldots,n\} \; , \; a\!\in\!\cA_i)$ is a 
vector of nonzero real numbers, and $(\ell_1,\ldots,\ell_n)$ is a mixed    
$n$-tuple of lifting functions for $\cA_1,\ldots,\cA_n$. Let 
$\Sigma_{\ell_1,\ldots,\ell_n}$ denote the 
resulting polyhedral subdivision of $\conv(\cA_1)+\cdots+\conv(\cA_n)$  
(as in Definition \ref{dfn:mixed}) and let 
$s_i\!:=\!(\sign(c_{i,a})\; | \; a\!\in\!\cA_i)$ for all $i$.  
Then, for all $t\!>\!0$ sufficiently small, the system of 
polynomials $\left(\sum\limits_{a\in \cA_1} c_{1,a} t^{\ell_1(a)}x^a,\ldots,
\sum\limits_{a\in \cA_n} c_{n,a} t^{\ell_n(a)}x^a\right)$  
has exactly $N$ roots in $\Rn_+$, where $N$ is the number of 
alternating cells of $(\Sigma_{\ell_1,\ldots,\ell_n},s_1,\ldots,s_n)$. \qed  
\end{sturmf} 

A final tool we will need is the {\em non-Archimedean Newton polytope},  
along with a recent refinement incorporating generalized phase. In particular, 
the definition and theorem below are special cases of a non-Archimedean 
analogue (see \cite{ai2}) of Sturmfel's result above.   
\begin{dfn} Given any complete non-Archimedean field $K$ with uniformizing 
parameter $\rho$, and any Laurent polynomial 
$f(x)\!:=\!\sum^m_{i=1}c_ix^{a_i}\!\in\!K[x^{\pm 1}_1,
\ldots,x^{\pm 1}_n]$, we define its {\em Newton polytope over $K$} to be 
$\newt_K(f)\!:=\!\conv\!\left\{(a_i,\ord\; c_i)\; | \; i\!\in\!\{1,\ldots,
m\} \right\}$. Also, the polynomial associated to summing the terms of $f$ 
corresponding to points of the form $(a_i,\ord\; c_i)$ lying on a lower face 
of $\newt_K(f)$, and replacing each coefficient $c$ by its first digit 
$\phi(c)$, is called a {\em lower polynomial}. 
\dia 
\end{dfn}

\noindent 
A remarkable fact true over non-Archimedean algebraically closed 
fields, but false over $\C$, is that the norms of roots of polynomials can be 
determined completely combinatorially: see Section \ref{sub:tropical} below and 
\cite{ekl}. What is less well-known is that, under 
certain conditions, the generalized phases can also be found by 
simply solving some lower binomial systems. Henceforth, 
we abuse notation slightly by setting $\ord(y_1,\ldots,y_n)\!:=\!
(\ord \; y_1,\ldots,\ord \; y_n)$. 
\begin{thm}
\label{thm:ai}
(Special case of \cite[Thm.\ 3.10 \& Prop.\ 4.4]{ai2}.) 
Suppose $K$ is a complete non-Archimedean field with residue field $\fk$ and  
uniformizer $\rho$. Also let 
$f_1,\ldots,f_n\!\in\!K[x^{\pm 1}_1,\ldots,x^{\pm 1}_n]$, $\hQ\!:=\!
\sum^n_{i=1}\newt_K(f_i)$, and let $(v,1)$ be an 
inner normal to a mixed lower\linebreak 
facet of $\hQ$ of the form $\hE\!:=\!\hE_1+\cdots+\hE_n$ 
where $\hE_i$ is a lower edge of $\newt_K(f_i)$ for all $i$.\linebreak 
Suppose also that the lower polynomials $g_1,\ldots,g_n$ corresponding to 
the normal $(v,1)$ are\linebreak 
\scalebox{.935}[1]{all binomials, and that $\pi\!\left(\hE\right)$ has 
standard Euclidean volume $1$. Then $F\!:=\!(f_1,\ldots,f_n)$ has 
$1$ or $0$}\linebreak 
\scalebox{.95}[1]{roots $\zeta\!\in\!(K^*)^n$ with $\ord\; \zeta\!=\!v$ and 
generalized phase $\theta\!\in\!(\fk^*)^n$ according as 
$g_1(\theta)\!=\cdots=\!g_n(\theta)\!=\!0$}\linebreak 
or not. In particular, $F$ has at 
most one root with valuation vector $v$. \qed  
\end{thm}

\noindent 
Note that while the number of roots with given $n$-tuple of first digits may 
depend on the uniformizer $\rho$ (see Proposition \ref{prop:iou} in Section 
\ref{sec:wrap}), the total number of roots with $\ord\; \zeta\!=\!v$ is 
independent of $\rho$. 
\begin{ex}
\label{ex:3d} 
Let $p$ be any prime, $n\!=\!3$, and let $(\cA_1,\cA_2,\cA_3)$ be the 
triple of supports for the system $G_p$ (see Theorem \ref{thm:family}). Also 
let $\ell_1,\ell_2,\ell_3$ be the respective liftings obtained by using the 
$p$-adic valuations of the coefficients of $G_p$. Lemma \ref{lemma:tri} then 
tells us that we obtain exactly $4$ mixed cells (two views of which 
are shown below), with corresponding 
lower facet normals $(1,0,0,1),(0,0,0,1),(-1,-1,0,1),(-2,-2,-1,1)$. In 
particular, the corresponding lower binomial systems are the following:\\  
\epsfig{file=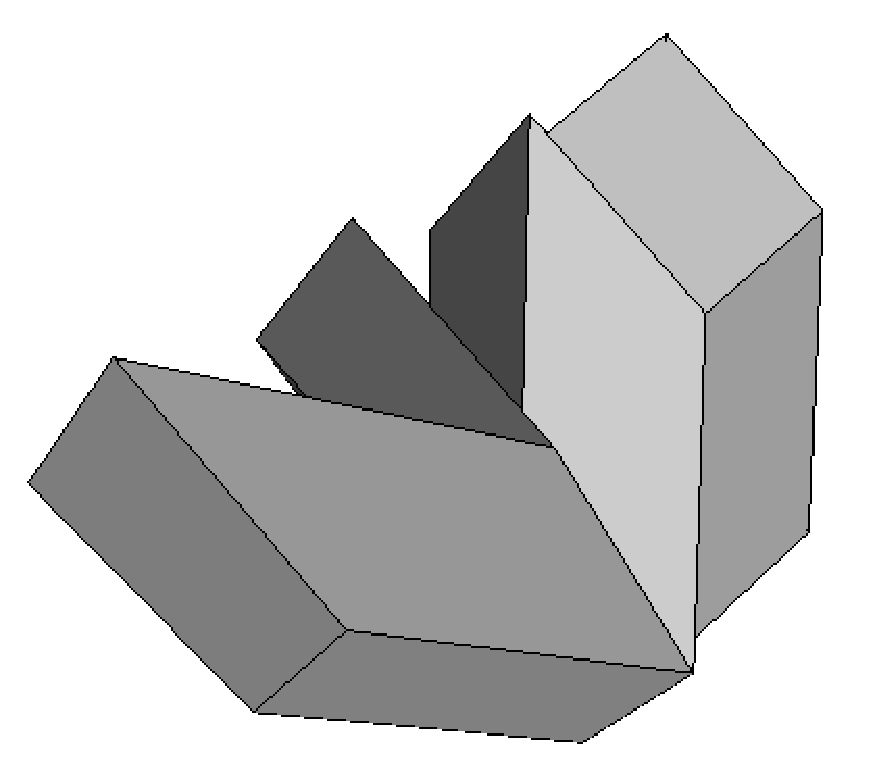,height=1in} 
\raisebox{1.2cm}{$\left. 
\begin{matrix} 
x_1x_2-1\\
x_2x_3-1\\
x_3-1
\end{matrix} \  \right| \left. \ 
\begin{matrix} 
x_1x_2-x^2_1\\
x_2x_3-1\\
x_3-1
\end{matrix} \ \right| \left. \ 
\begin{matrix} 
x_1x_2-x^2_1\\
x_2x_3-x^2_1\\
x_3-1
\end{matrix} \ \right| \ 
\begin{matrix} 
x_1x_2-x^2_1\\
x_2x_3-x^2_1\\
x_3-x^2_1
\end{matrix}$} 
\epsfig{file=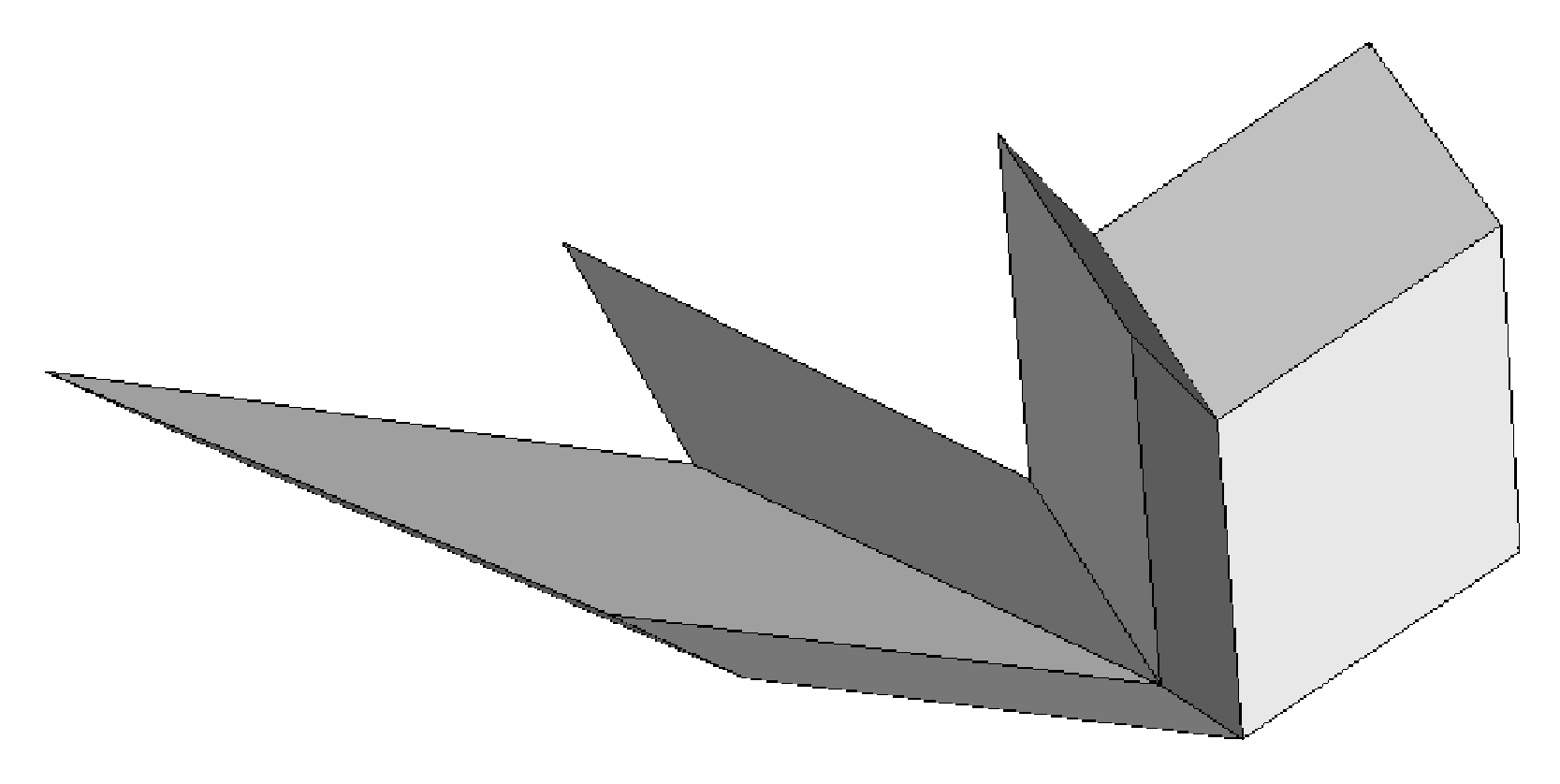,height=1in}\linebreak
Each mixed cell has volume $1$, and each corresponding binomial system has 
unique solution $(1,1,1)\!\in\!(\F^*_p)^3$. 
Theorem \ref{thm:ai} then tells us that the roots of $G_p$ in $(\Q^*_p)^3$ are 
of the following form: \ \ $(p(1+O(p)),1+O(p),1+O(p))$, \  \   
$(1+O(p),1+O(p),1+O(p))$,\linebreak 
\scalebox{.96}[1]{$\left(p^{-1}(1+O(p)),p^{-1}(1+O(p)),1+O(p)\right)$, and  
$\left(p^{-2}(1+O(p)),p^{-2}(1+O(p)),p^{-1}(1+O(p))\right)$. \dia}   
\end{ex} 

\subsection{Some Tropical Visualizations} 
\label{sub:tropical} 
A beautiful theorem of Kapranov tells us that, for 
non-Archimedean $K$, we can use polyhedral combinatorics to 
efficiently compute the valuations of the roots of any polynomial. 
\begin{dfn}
\label{dfn:amoeba} 
For any complete algebraically closed field $K$ 
and $f\!\in\!K\!\left[x^{\pm 1}_1,\ldots,x^{\pm 1}_n\right]$   
we set $Z^*_K(f)\!:=\!\{x\!\in\!(K^*)^n\; | \; 
f(x)\!=\!0\}$. Also, for any subset $S\!\subseteq\!\Rn$, we let $\bar{S}$ 
denote the closure of $S$ in the Euclidean topology. 
Finally, if $K$ is also non-Archimedean, then  
we define the {\em tropical variety of $f$ over $K$}, $\trop_K(f)$, to be the 
closure in $\Rn$ of\\ 
\mbox{}\hfill $\{(v_1,\ldots,v_n)\!\in\!\Rn\; | \; 
(v_1,\ldots,v_n,1) \text{ is an inner edge normal of } \newt_K(f) 
\}$\hfill\mbox{}\dia 
\end{dfn} 
\begin{rem} 
$\trop_K(f)$ is sometimes equivalently defined in terms of max-plus  
semi-rings (see, e.g., \cite{macsturmf}). \dia  
\end{rem} 
\begin{kapra} \cite{ekl} 
For any complete,\linebreak non-Archimedean algebraically closed field 
$K$, we have   
$\overline{\ord\!\left(Z^*_K(f)\right)}\!=\!\trop_K(f)$. \qed 
\end{kapra}

We now illustrate these ideas through our earlier examples. 
Returning to Example \ref{ex:g2}, 
the underlying tropical varieties (or closures of $\ord\!\left(
Z^*_L(g_1)\right)$ and $\ord\!\left(Z^*_L(g_2)\right)$ for 
$L\!\in\left\{\overline{\Q}_p,\overline{\F_q((t))}\right\}$) intersect in 
exactly $3$ points as illustrated below, on the left. (The tropical varieties 
for the first and second polynomials are respectively colored in solid red 
and dashed blue.) The right-hand illustration below shows the corresponding 
plots when $L\!=\!\C$ and $\eps\!=\!1/4$, with their intersection darkened 
slightly.\\ 
\mbox{}\hfill \epsfig{file=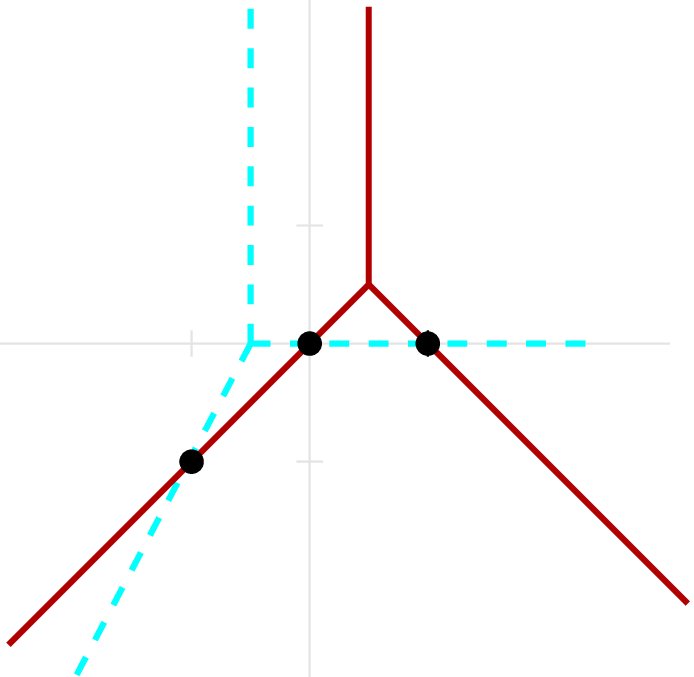,height=1.8in,clip=}\hspace{1.5cm} 
\epsfig{file=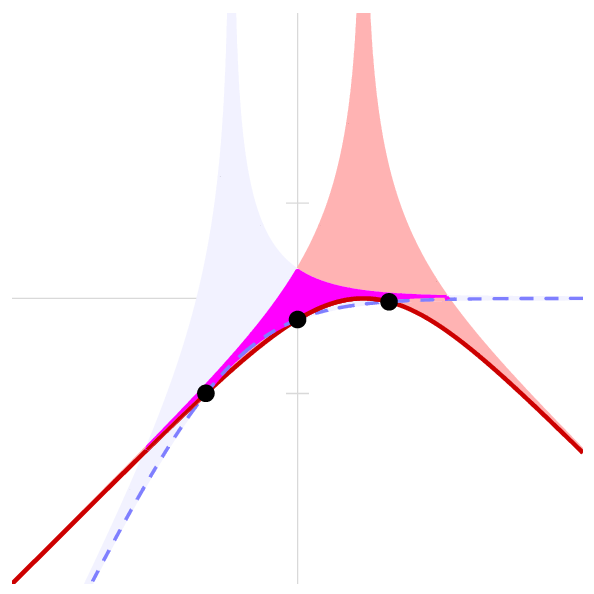,height=1.8in,clip=}\hfill\mbox{} \\ 
Note that the images of the corresponding positive 
zero sets under the (complex) $\ord$ map are drawn as even darker 
curves (with $3$ marked intersections) in the right-hand illustration above. 
The negative of the image of a complex algebraic set under the complex $\ord$ 
map is usually called an {\em amoeba} \cite{passare}.  

Returning to Example \ref{ex:3d}, the resulting tropical varieties are 
illustrated below (without translucency on the left, with translucency on the 
right): \\ 
\mbox{}\hfill 
\epsfig{file=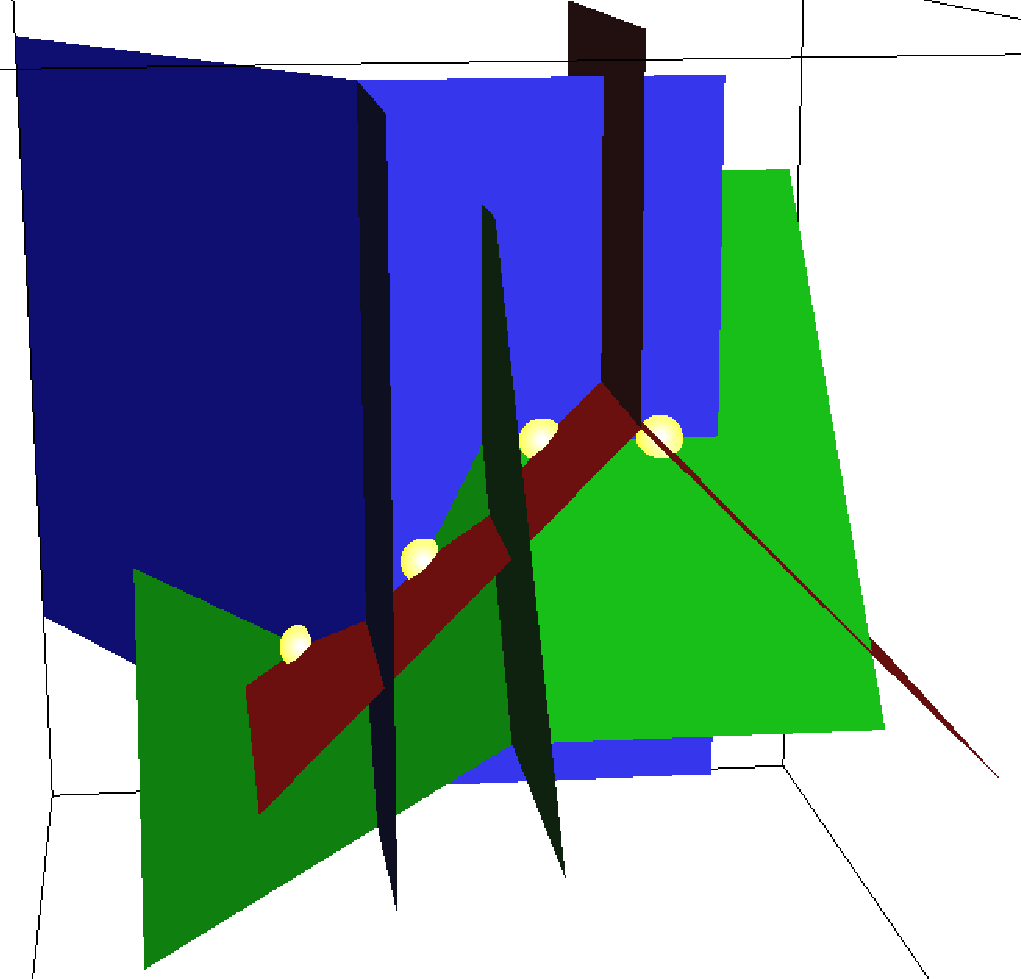,height=2.8in,clip=} \hfill 
\epsfig{file=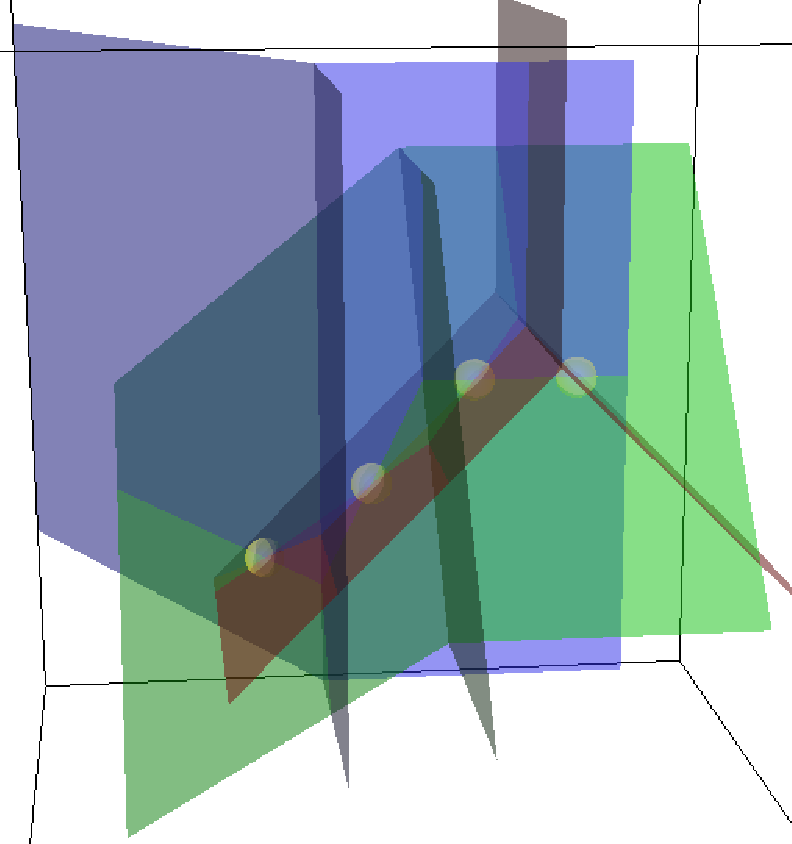,height=2.8in,clip=}
\hfill \mbox{} \\ 
Note that each tropical variety above is a polyhedral complex of 
codimension $1$, and that all the top-dimensional faces are unbounded,  
even though they are truncated in the illustrations. 

\newpage 
\section{Proving our Main Results}  
\label{sec:proof} 
\subsection{Theorem \ref{thm:big}: The Universal Lower Bound} 
\label{sub:big} \mbox{}\\ 
First note that since $Y_L(n,k)$ is integer-valued when finite, 
$Y_L(n,k)$ is actually attained by some $(n+k)$-nomial $n\times n$ system over 
$L$ when $Y_L(n,k)$ is finite. 
 
Now, any $n\times n$ polynomial system of the form 
$(b(x_1),\ldots,b(x_{n-1}),r(x_n))$ --- with $b\!\in\!L[x_1]$ a binomial and 
$r\!\in\!L[x_1]$ a trinomial, both possessing nonzero constant terms --- is 
clearly an $(n+2)$-nomial $n\times n$ system. 
So we immediately obtain $Y_L(n,2)\!\geq\!Y_L(1,2)Y_L(1,1)^{n-1}$ 
simply by picking $b$ and $r$ (via Theorem \ref{thm:sum} and
Remark \ref{rem:low}) to have maximally many roots over $L$ with all 
coordinates of generalized phase $1$. That $Y_L(n,2)\!\geq\!n+1$ follows 
immediately from Theorem \ref{thm:family}, so we obtain the first asserted 
inequality. 

The remaining lower bounds for $Y_L(n,k)$ follow from similar  
concatenation tricks. First, note that 
any $n\times n$ polynomial system of the form 
$(b(x_1),\ldots,b(x_{n-k+1}),r(x_{n-k+2}),\ldots,r(x_n))$ 
is clearly an $(n+k)$-nomial $n\times n$ system. So, specializing 
$b$ and $r$ appropriately once again, the inequality 
$Y_L(n,k)\!\geq\!Y_L(1,1)^{n-k+1}Y_L(1,2)^{k-1}$ 
holds for $n\!\geq\!k-1$. 

A slightly more intricate construction gives our next lower bound:  
letting $F_n(x_1,\ldots,x_n)$ denote an $(n+2)$-nomial $n\times n$ system 
over $L$ possessing a nonzero constant term, observe that 
when $k-1\!\leq\!n$ and $\ell\!:=\!\lfloor \frac{n}{k-1}\rfloor$, the 
block-diagonal system $F$ defined by\\
\mbox{}\hfill  
$F_\ell(x_{1,1},\ldots,x_{1,\ell}),\ldots,F_\ell(x_{k-1-[n]_{k-1},1},
\ldots,x_{k-1-[n]_{k-1},\ell})$,\hfill\mbox{}\\ 
\mbox{}\hfill $F_{\ell+1}(y_{1,1},\ldots,y_{1,\ell+1}),
\ldots,F_{\ell+1}(y_{[n]_{k-1},1},\ldots,y_{[n]_{k-1},\ell+1})$ \hfill
\mbox{}\\ 
involves exactly $(k-1-[n]_{k-1})\ell+[n]_{k-1}(\ell+1)\!=\!(k-1)
\ell+[n]_{k-1}\!=\!n$ 
variables, and $n$ polynomials via the same calculation. Also, the total 
number of distinct exponent vectors of $F$ is exactly\\ 
\mbox{}\hfill $(k-1-[n]_{k-1})(\ell+2)+[n]_{k-1}(\ell+3)-(k-1)+1\!=\!
(k-1)\ell+[n]_{k-1}+2(k-1)-k+2\!=\!n+k$,\hfill\mbox{}\\  
since all the polynomials share a nonzero constant term. 
Furthermore, any ordered $n$-tuple consisting of $k-1-[n]_{k-1}$ 
non-degenerate roots of $F_\ell$ in $L^\ell$ 
followed by $[n]_{k-1}$ non-degenerate 
roots of $F_{\ell+1}$ in $L^{\ell+1}$ (with all coordinates having 
generalized phase $1$) is clearly a non-degenerate root of $F$ 
in $L^n$ with all coordinates having generalized phase $1$. Picking 
$F_\ell$ and $F_{\ell+1}$ to be appropriate specializations of the 
systems from Theorem \ref{thm:family}, we thus obtain 
$Y_L(n,k)\!\geq\!Y_L\!\left(\left\lfloor \frac{n}{k-1}\right\rfloor,2
\right)^{k-1-[n]_{k-1}}Y_L\!\left(\left\lfloor \frac{n}{k-1}\right
\rfloor+1,2\right)^{[n]_{k-1}}$. So the case $n\!\geq\!k-1$ is done. 

Now simply note that any $n\times n$ system of the form\\ 
\mbox{}\hfill  
$(m(x_1),\ldots,m(x_{n-[k-1]_n}),\mu(y_1),\ldots,\mu(y_{[k-1]_n}))$\hfill
\mbox{}\\ --- with $m\!\in\!L[x_1]$ an $\ell$-nomial, $\mu\!\in\!L[y_1]$ 
an $(\ell+1)$-nomial, $\ell\!:=\!\lfloor \frac{n+k-1}{n}\rfloor$, and\linebreak 
$n\!\leq\!k-1$ ---  is easily verified to be an $(n+k)$-nomial $n\times n$ 
system. So picking $m$ and $\mu$ to have maximally many roots with generalized 
phase $1$, we immediately obtain $Y_L(n,k)\!\geq\!Y_L\!\left(1,\left\lfloor 
\frac{n+k-1}{n}\right\rfloor\right)^{n-[k-1]_{n}}Y_L\!\left(1,\left\lfloor 
\frac{n+k-1}{n} \right\rfloor+1\right)^{[k-1]_{n}}$ for $n\!\leq\!k-1$. 

To conclude, the entries in our table are simply specializations of 
our recursive lower bounds using the explicit values given by 
Theorem \ref{thm:sum}. \qed 

\subsection{Theorem \ref{thm:family}: Fewnomials Systems with Many Roots 
Universally} 
\label{sub:family} \mbox{}\\ 
First note that all the roots of $G_\eps$ in $\bar{L}^n$ lie in 
$\left(\bar{L}^*\right)^n$. (Clearly, setting any $x_i\!=\!0$ results in a 
pair of univariate polynomials having no roots in common, or a 
nonzero constant being equal to zero.) 
Let $(g_1,\ldots,g_n)\!:=\!G_\eps$ and let $\cA$ denote the matrix whose 
columns are the vectors in the union of the supports of the $g_i$. More 
precisely, $\cA$ is the $n\times (n+2)$ matrix below:\\ 

\vspace{-.8cm}
\noindent
\begin{minipage}[t]{.2\textwidth}
\vspace{0pt}
\scalebox{.8}[.48]{$\begin{bmatrix} 
0      & 2 & 1 & 0 &        &   & \\ 
       &   & 1 & 1 &        &   & \\ 
       &   &   & 1 &        &   & \\
       &   &   &   & \ddots &   & \\
       &   &   &   &        & 1 &  \\
       &   &   &   &        & 1 & 1
\end{bmatrix}$} 
\end{minipage} \hspace{.5cm}
\begin{minipage}[t]{.75\textwidth}
\vspace{0pt}
Now let $\bA$ denote the $(n+1)\times (n+2)$ matrix obtained by 
appending a row of $1$s to the top of $\cA$. It is then 
easily checked that $\bA$ has right null-space of dimension $1$, 
generated by the transpose of 
$b\!:=\!(b_1,\ldots,b_{n+2})\:=\!(-1,(-1)^n,(-1)^{n+1} 2,
\ldots,(-1)^{n+n}2)$. 
Let us rewrite\linebreak 
the equation $g_i\!=\!0$ as $x^{a_{i+2}}\!=\!\beta_i(x^2_1)$,
where $a_i$ denotes the $i\thth$ column\linebreak of $\cA$ and
$\beta_i$ is a suitable degree one polynomial with coefficients that\linebreak 
\end{minipage}  

\vspace{-.37cm} 
\noindent 
are powers of $\eps$. Since the entries of $b$ sum to $0$, we then easily
obtain that\\ 
\mbox{}\hfill $1^{b_1}u^{b_2}\beta_1(u)^{b_3}\cdots\beta_n(u)^{b_{n+2}}\!=\!1$
\hfill\mbox{}\\ 
when $\zeta\!=\!(\zeta_1,\ldots,\zeta_n)$ is a root of 
$G_\eps$ in $\left(\bL^*\right)^n$ and $u\!:=\!\zeta^2_1$. 
In other words, the degree $n+1$ polynomial 
$R_n(u)$ from Lemma \ref{lemma:kait} must vanish. 
Furthermore, the value of $\zeta_n$ is uniquely determined by the value of $u$, 
thanks to the equation $g_n\!=\!0$. Proceeding with the remaining 
equations $g_{n-1}\!=\!0,\ldots,g_1\!=\!0$ we see that the same holds for  
$\zeta_{n-1},\ldots,\zeta_2$ and $\zeta_1$ successively. So $G_\eps$ 
has no more than $n+1$ roots, counting multiplicities, in 
$\left(\bL^*\right)^n$. 
Note in particular that by Lemma \ref{lemma:tri}, combined with 
{\em Bernstein's Theorem} (over a general algebraically closed field 
 \cite{bernie,danilov}), $G_\eps$ having 
at least $n+1$ distinct roots in $\left(\bL^*\right)^n$ implies that there are 
{\em exactly} $n+1$ roots in $\left(\bL^*\right)^n$ and they are all 
non-degenerate. 

To finally prove the first part of our theorem, we separate the Archimedean 
and non-Archimedean cases: when $L\!=\!\R$ we immediately obtain,  
from Lemma \ref{lemma:tri} and Sturmfels' Theorem, that $G_\eps$ has at 
least $n+1$ positive roots for $\eps\!>\!0$ sufficiently small. (This 
trivially implies the $L\!=\!\C$ case as well.) 

For the non-Archimedean case, Lemma \ref{lemma:tri} and Theorem \ref{thm:ai} 
immediately imply that, when $\phi(\eps)\!=\!1$ and 
$\ord \; \eps\!\geq\!1$, $G_\eps$ has at least $n+1$ roots 
in $L^n$ with all coordinates having generalized phase $1$. In particular,
for each vector $v_j$ from Lemma \ref{lemma:tri},  
it is easily checked that 
$(1,\ldots,1)$ is a root of the corresponding lower binomial system of 
$G_\eps$ over the residue field of $L$.

The only assertion left to prove is that $G_{1/4}$ has exactly $n+1$ roots in 
the positive orthant, and this follows from Lemma \ref{lemma:kait}. \qed   

\subsection{Proof of Lemma \ref{lemma:kait}} 
\label{sub:kait} 
Let us first define $A_n$ and $B_n$ respectively as\\  
$u (1+\eps u)^2 (1+\eps^5 u)^2\cdots 
(1+\eps^{4\lfloor n/2\rfloor-3}u)^2$ and 
$(\eps+u)^2 (1+\eps^3 u)^2 (1+\eps^7 u)^2 \cdots 
 (1+\eps^{4\lceil n/2\rceil-5}u)^2$. 
Clearly, $R_n\!=\!A_n-B_n$. 
\begin{lemma}\label{lemma:k1} \scalebox{.95}[1]{Assume $\eps\!=\!1/4$.
Then, for all $n\!\geq\!2$, we have
$\displaystyle{R_n\!\left(16^{n-2}/u\right)\!=\!\left(\frac{-4^{n-2}}{u}
\right)^{n+1}R_n(u)}$.} 

\vspace{-.2cm}
\noindent
Also, for all even $n\!\geq\!2$, we have
$R_n(4^{n-2})=0$.
\end{lemma}
\begin{lemma}\label{lemma:k2} Assume $\eps\!=\!1/4$ and 
consider $R_n$ as a function on $\R$. Then,
for all $n\geq2$, we have (a) $R_n(0)\!<\!0$ and
(b) $(-1)^\ell R_n(16^\ell/4)\!>\!0$ for all
$\ell\!\in\!\{0,\ldots,\ceil{n/2}-1\}$.
\end{lemma}

\noindent
These subsidiary lemmata are proved in Section \ref{sec:wrap} below.

Returning to the proof of Lemma \ref{lemma:kait},  
we now consider two exclusive cases. 

\medskip
\noindent
{\bf Real Case:} 
By Lemma \ref{lemma:k2}, $R_n$ has
$\ceil{\frac{n}{2}}-1$ sign changes in the open interval
$\left(0,\frac{16^{\ceil{n/2}-1}}{4}\right)$. So by the Intermediate
Value Theorem, $R_n$ has
$\ceil{\frac{n}{2}}-1$ roots in this interval. By Lemma \ref{lemma:k1}, for
every such root $\zeta$, $\frac{16^{n-2}}{\zeta}$ yields a new root. When
$n$ is odd, this gives us $2(\ceil{\frac{n}{2}}-1)=n+1$ positive roots.  When
$n$ is even, we get $n$ positive roots and, by Lemma \ref{lemma:k1},
the new positive root $4^{n-2}$.  So $R_n$ has $n+1$ positive
roots. \qed

\medskip 
\noindent  
{\bf Non-Archimedean Case:} 

\vspace{-.3cm} 
\noindent 
\begin{minipage}[t]{.18\textwidth}
\vspace{0pt}
\epsfig{file=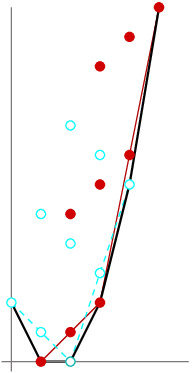,height=2.2in}
\end{minipage} 
\begin{minipage}[t]{.82\textwidth}
\vspace{0pt} 
\mbox{}\hspace{.5cm}
While this case is already implicit 
in the proof of Theorem \ref{thm:family}, one can form a direct argument 
starting from Newton polygons: For $L\!\in\!\{\Q_p,\F_q((t))\}$ (and thus 
$\eps\!\in\!\{p,t\}$ respectively), we easily obtain that 
$P\!:=\!\newt_L(A_n)$ has exactly $1+\lfloor n/2\rfloor$ lower edges, 
$Q\!:=\!\newt_L(B_n)$ has exactly $\lceil n/2\rceil$ lower edges, and the 
vertices of $P$ and $Q$ interlace. 
(The supports of $A_4$ and $B_4$ are drawn, respectively as red (filled) and 
blue (unfilled) circles, at left.) More precisely, 
$\newt_L(R_n)\!=\!\conv(P\cup Q)$ has exactly $n+1$ 
lower edges, each having horizontal length $1$. In particular, 
$\{(1,1),(0,1),\ldots,(1-n,1)\}$ is a representative set of inner normals 
for the lower edges, and each corresponding 
lower binomial is a degree one polynomial with  
pair of coefficients $(\pm 1,\mp 1)$. Also, for any 
$i\!\in\!\{1,0,\ldots,1-n\}$,\linebreak  
\end{minipage}  

\vspace{-.35cm} 
\noindent 
we can find a $d_i\!\in\!\Z$ such that $\eps^{d_i}R_n(\eps^i u)\!=\!\pm 1 \mp u 
+ O(\eps)$. So by Hensel's 
Lemma, $R_n$ has exactly $n+1$ roots in $\Q_p$ (resp.\ $\F_p((t))$) when 
$\eps\!=\!p$ (resp.\ $\eps\!=\!t$), and each such root has first digit $1$. 
\qed 

\subsection{Proof of Lemma \ref{lemma:tri}} 
\label{sub:tri} 
By Theorem \ref{thm:mixed} our mixed volume in question is bounded 
above by $n!\vol(Q)$ where $Q$ is the polytope with vertices the 
columns of the matrix $\cA$ from the proof of Theorem \ref{thm:family}. 
The vertices of $Q$ form a {\em circuit}, and the signs of the 
entries of the vector $b$ from the proof of Theorem \ref{thm:mixed}  
thereby encode an explicit triangulation of $Q$ 
(see, e.g., \cite[Prop.\ 1.2, pg.\ 217]{gkz94}). More precisely, defining 
$Q(i)$ to be the convex hull of the points corresponding to all the columns 
of $\cA$ {\em except} for the $i\thth$ column, we obtain that 
$\left\{Q(2),Q(4),\ldots,Q\!\left(2\left\lfloor\frac{n+2}{2}\right\rfloor
\right)\right\}$ (for 
$n$ even) and $\left\{Q(3),Q(5),\ldots,Q\!\left(2\left\lceil\frac{n+2}{2}
\right\rceil-1\right)\right\}$ (for $n$ odd) form the simplices of a 
triangulation of $Q$. Note in particular that the volume of $Q(i)$ is exactly 
$1/n!$ times the absolute value of the determinant of the submatrix of 
$\cA$ obtained by deleting the first and $i\thth$ columns. Note also that 
this submatrix is block-diagonal with exactly $2$ blocks: an 
$(i-2)\times (i-2)$ upper-left upper-triangular block and an 
$(n-i+2)\times (n-i+2)$ lower-right lower-triangular block. It is then 
clear that $\vol(Q(i))$ is $1$ or $2$, according as $i\!=\!2$ or  
$i\!\geq\!3$. So $\vol(Q)$ is then $1+2\left(\left\lfloor\frac{n+2}{2}
\right\rfloor-1\right)\!=\!n+1$ (when $n$ is even) or 
$2\left(\left\lceil\frac{n+2}{2}\right\rceil-1\right)\!=\!n+1$ (when 
$n$ is odd). 

Since any $n$-tuple of columns chosen from the 
last $n+1$ columns of $\cA$ is linearly independent, 
each cell $\pi\!\left(\hP_j\right)$ has positive volume. 
(The linear independence follows directly from our preceding block diagonal 
characterization of certain submatrices of $\cA$.) So once we show that each 
such cell is distinct, we immediately obtain that our mixed volume is at least 
$n+1$ and thus equal to $n+1$. Toward this end, we now check that each $v_j$ 
is indeed an inner normal to $\hP_j$.  

For any $i\!\in\!\{1,\ldots,n\}$ let $\hA_i\!=\!(\alpha_i,\beta_i,\gamma_i)$ 
denote the triple of vertices 
of the triangle $\hT_i$, ordered so that $\pi(\alpha_i)\!=\!\bO$ and 
$\pi(\beta_i)\!=\!2e_1$. It then clearly suffices to prove that, for any 
$j\!\in\!\{0,\ldots,n\}$, the inner product $v_j\cdot x$ is minimized on each 
$\hA_i$ exactly at the vertices of the edge $\hE_{i,s}$, where $s$ is $1$ or 
$0$ according as $i\!\leq\!j$ or $i\!\geq\!j+1$. Equivalently, this means 
that the minimum values in the triple 
$(v_j\cdot\alpha_i,v_j\cdot\beta_i,v_j\cdot
\gamma_i)$ must occur exactly at the second and third (resp.\ first and third) 
coordinates when $i\!\leq\!j$ (resp.\ $i\!\geq\!j+1$). This follows from a 
direct but tedious computation that we omit. \qed

\subsection{Proofs for Example \ref{ex:slp}} 
\label{sub:slp} 
The assertion on the degree of $\frac{h_{n}(x_1)}{x_1(1-x_1)}$ is 
obvious from the recurrence for $h_{n}$. 
The upper bound on $\tau\!\left(\frac{h_{n}(x_1)}{x_1(1-x_1)}\right)$ 
follows easily from 
recursive squaring. 
In particular, since $\tau(c_S)\!\leq\!2\log_2 c_S$,  
we easily obtain $\tau(c_S)\!=\!O(\#S\log k)$. Expressing 
$c^{3^{n-1}}_S\!=\!(\cdots (c^3_S)^3 \cdots)^3$, it is then clear that 
$\tau\!\left(c^{3^{n-1}}_S\right)\!=\!O(n+\#S\log k)$. Observing that we can 
easily evaluate $\frac{h_{n}(x_1)}
{x_1(1-x_1)}$ by simply replacing $h_{2}$ by $c_S-h_{1}$ in the recurrence 
for $h_{n}$, we arrive at our bound for $\tau\!\left(\frac{h_{n}(x_1)}
{x_1(1-x_1)}\right)$. Note also that by construction, 
$\frac{h_{n}(x_1)}{x_1(1-x_1)}$ does not vanish at $0$ or $1$, but 
does vanish at every other root of $h_{n}$. 

We now focus on counting the roots of $\frac{h_{n}(x_1)}{x_1(1-x_1)}$ 
in the rings $\Z_p$ for $p\!\in\!S$. From our last 
observations, it clearly suffices to show that, for all $n\!\geq\!1$, 
$h_{n}$ has exactly $2^n$ roots in $\Z_p$ for each 
$p\!\in\!S$. We do this by induction, using the following 
refined induction hypothesis:  

\smallskip 
\noindent 
\mbox{}\hspace{1cm}For any prime $p\!\in\!S$, 
$h_{n}$ has exactly $2^n$ distinct roots in $\Z_p$. Furthermore, these roots\\ 
\mbox{}\hspace{1cm}are distinct mod $p^{3^{n-1}}$ 
and, for any such root $\zeta$, we have 
$\ord \; h'_n(\zeta)\!=\!\frac{3^{n-1}-1}{2}$.\\   
The case $n\!=\!1$ is clear. One also observes $h'_{1}(x_1)
\!=\!1-2x_1$, and $h'_{n+1} \!=\!(c^{3^{n-1}}_S-h_{n})h'_{n}$ for 
all $n\!\geq\!1$. So let us now assume the 
induction hypothesis for any particular $n$ and prove the case $n+1$. In 
particular, let $\zeta\!\in\!\Z_p$ be any of the $2^n$ roots of $h_{n}$. 
The derivatives of $h_{n}$ and $c^{3^{n-1}}_S-h_{n}$ differ 
only by sign mod $p^{3^{n-1}}$, so by Hensel's Lemma (combined with our 
induction hypothesis), $c^{3^{n-1}}_S-h_{n}$ also has $2^n$ distinct 
roots in $\Z_p$. However, the roots of $c^{3^{n-1}}_S-h_{n}$ in $\Z_p$ are all 
distinct from the roots of $h_{n}$ in $\Z_p$: this is because 
$c^{3^{n-1}}_S-h_{n}$ is nonzero at every root of $h_{n}(x_1)$ mod 
$p^{3^{n-1}+1}$. So $h_{n+1}$ then clearly has $2^{n+1}$ distinct roots in 
$\Z_p$, and these roots remain distinct mod $p^{3^n}$. Furthermore, by our 
recurrence for $h'_{n}$, the $p$-adic 
valuation of $h'_{n+1}$ is exactly $3^{n-1}+\frac{3^{n-1}-1}{2}\!=\!
\frac{3^n-1}{2}$. So our induction is complete. 

To see that $\frac{h_{n}(x_1)}{x_1(1-x_1)}$ has no real roots, first note 
that $x_1(1-x_1)$ is strictly increasing on $(-\infty,1/2)$, 
strictly decreasing on $(1/2,+\infty)$, and attains a unique 
maximum of $1/4$ at $x_1\!=\!1/2$. Since 
$c_S\!\geq\!2$, we also clearly obtain that $c_S-x_1(1-x_1)$ has range 
contained in $[3/4,+\infty)$, with minimum occuring at $x_1\!=\!1/2$. 
More generally, our recurrence for $h'_{n}$ implies that any critical 
point $\zeta\!\in\!\R$ of $h'_n$, other than a critical point of $h_{n-1}$, 
must satisfy $c^{3^{n-1}}_S\!=\!2h_{n-1}(\zeta)$. So, in particular, 
$h_{2}$ has the same regions of strict increase and strict decrease as 
$h_{1}$, and thus $h_{2}$ has maximum $\leq\!3/8$. Proceeding by 
induction, we see thus see that $h_{n}$ has no critical points other than 
$1/2$ and thus no real roots other than $0$ and $1$. Moreover, the latter 
roots occur with multiplicity $1$ from the obvious recursive factorization of 
$h_{n}$. So $\frac{h_{n}(x_1)}{x_1(1-x_1)}$ has no real roots. \qed  

\section{Wrapping up: Invariance of $Y_L(n,k)$, and the Proofs of 
Proposition \ref{prop:simp},  
Theorem \ref{thm:hah}, and Lemmata \ref{lemma:k1} and \ref{lemma:k2}}  
\label{sec:wrap} 
Let us now see how the value of $Y_L(n,k)$ depends weakly (if at all) on the 
underlying uniformizer, and how counting roots with coordinates of generalized 
phase $1$ is as good as counting roots in any other direction. In what 
follows, we let $W_L(n,k)$ denote the supremum, over all $(n+k)$-nomial 
$n\times n$ systems $F$ over $L$, of the {\em total} number of non-degenerate 
roots of $F$ in $(L^*)^n$. 
\begin{prop}
\label{prop:iou}
\mbox{}\\ 
(1) For $L$ any finite extension of $\Q_p$, and $n,k\!\geq\!1$, the value 
of $Y_L(n,k)$ in Definition \ref{dfn:omega} is\linebreak 
\mbox{}\hspace{.7cm}independent of the choice of 
uniformizer $\rho$. Also, the same holds for $L\!=\!\F_q((t))$\linebreak 
\mbox{}\hspace{.7cm}when $n\!=\!1$. \\ 
(2) $Y_L(n,k)$ counts the supremum of the number of roots in {\em any} fixed 
angular direction in the\linebreak 
\mbox{}\hspace{.7cm}following sense: let $\theta_1,\ldots,\theta_n$ be elements 
of the complex unit circle, elements of $\{\pm 1\}$, or\linebreak 
\mbox{}\hspace{.7cm}units in the residue 
field of $L$, according as $L$ is $\C$, $\R$, or non-Archimedean. Also, 
letting\linebreak
\mbox{}\hspace{.7cm}$F$ and $G$ denote $(n+k)$-nomial $n\times n$ systems over 
$L$, there is an $F$ with exactly $N$ non-\linebreak 
\mbox{}\hspace{.7cm}degenerate roots 
$(\zeta_1,\ldots,\zeta_n)\!\in\!L^n$ satisfying $\phi(\zeta_i)\!=\!\theta_i$ 
for all $i$ if and only if there is a $G$\linebreak 
\mbox{}\hspace{.7cm}\scalebox{.98}[1]{with exactly $N$ non-degenerate 
roots in $L^n$ with all coordinates having generalized phase $1$.}\\  
(3) $W_\C(n,k)\!=\!+\infty$, $W_\R(n,k)\!=\!2^nY_\R(n,k)$, and 
$W_L(n,k)\!=\!(q_L-1)^nY_L(n,k)$ for any 
finite\linebreak 
\mbox{}\hspace{.7cm}extension $L$ of $\Q_p$ with residue 
field cardinality $q_L$.  Also, we have\\ 
\mbox{}\hspace{2.7cm}
$W_{\F_q((t))}(n,k)\leq (q-1)^n Y_{\F_q((t))}(n,k)\leq (q-1)^n
W_{\F_q((t))}(n,k)$. 
\end{prop} 

\noindent 
{\bf Proof:} \\ 
{\bf Assertion (2):} To prove independence of direction, fix a 
uniformizer $\rho$ once and for all (for the non-Archimedean case) and assume 
$F$ has exactly $N$ non-degenerate roots $(\zeta_1,\ldots,\zeta_n)\!\in\!L^n$ 
satisfying $\phi(\zeta_i)\!=\!\theta_i$ for all $i$. Defining 
$G(x_1,\ldots,x_n)
\!=\!F(t_1x_1,\ldots,t_nx_n)$ for any $t_1,\ldots,t_n$ of valuation $0$
with $\phi(t_i)\!=\!\theta_i$ for all $i$, we then clearly obtain a
suitable $G$ with exactly $N$ non-degenerate roots with all coordinates
having generalized phase $1$. The preceding substitutions can also be inverted 
to give the converse direction, so we obtain independence of direction, and  
(in the non-Archimedean case) for any $\rho$. \qed  

\smallskip
\noindent
{\bf Assertion (3):} The first equality was already observed in Section 
\ref{sub:up}. 

Now recall that any $y\!\in\!\R^*$ (resp.\ $y\!\in\!L$, $y\!\in\!\F_q((t))$) 
can be written in the form $y\!=\!uz$ where 
$u\!\in\!\{\pm 1\}$ (resp.\ $u$ is a unit in the residue field of $L$ or 
$u\!\in\!\F^*_q$), $|y|\!=\!|z|$, and $z$ has generalized phase $1$. 
So Assertion (2) then 
immediately implies $W_\R(n,k)\!\leq\!2^nY_\R(n,k)$,\linebreak  
$W_L(n,k)\!\leq\!(q_L-1)^nY_L(n,k)$, 
and $W_{\F_q((t))}(n,k)\!\leq\!(q-1)^nY_{\F_q((t))}(n,k)$. Note also 
that $Y_{\F_q((t))}(n,k)\!\leq\!W_{\F_q((t))}(n,k)$, independent of the 
underlying uniformizer. 

So now we need only prove $W_\R(n,k)\!\geq\!2^nY_\R(n,k)$ 
and $W_L(n,k)\!\geq\!(q_L-1)^nY_L(n,k)$. Toward this end, note 
that for any $F$ with $N$ non-degenerate roots in $\R^n$ (resp.\ $L^n$), 
with all coordinates of generalized phase $1$, the  
substitution $x_i\!=\!y^2_i$ (resp.\ $x_i\!=\!y^{q_L}_i$) for all $i$  
yields a new system with exactly $N$ non-degenerate roots in $\R^n$ (resp.\ 
$L^n$) with $n$-tuple of generalized phases $(\theta_1,\ldots,\theta_n)$ 
for {\em any} $\theta_1,\ldots,\theta_n$ in $\{\pm 1\}$ (resp.\ units in the 
residue field). Clearly then, $W_\R(n,k)\!\geq\!2^nY_\R(n,k)$ and 
$W_L(n,k)\!\geq\!(q_L-1)^nY_L(n,k)$. \qed  

\smallskip 
\noindent 
{\bf Assertion (1):} For $L$ as in the first part, Assertion (3) tells us that 
$Y_L(n,k)\!=\!\frac{W_L(n,k)}{(q_L-1)^n}$ where 
$q_L$ is the residue field cardinality of $L$. $W_L(n,k)$ is independent 
of $\rho$, so the first part is proved. The second assertion follows 
immediately from Section 2 of \cite{poonen}. \qed 

\subsection{Proof of Proposition \ref{prop:simp}}   
First note that by Gaussian elimination, $k\!\leq\!0$ immediately 
implies that any $(n+k)$-nomial $n\times n$ system is either equivalent to 
an $n\times n$ system where all the polynomials are monomials or an 
$n\times n$ system with at least one polynomial identically zero. Neither 
type of system can have a root in  
$(L^*)^n$ with Jacobian of rank $n$. So we obtain the first equality. 

Similarly, any $(n+1)$-nomial $n\times n$ system is either equivalent to 
an $n\times n$ system consisting solely of binomials or an 
$n\times n$ system with at least polynomial having $1$ or fewer monomial 
terms. The latter type of system can not have a root in
$(L^*)^n$ with Jacobian of rank $n$, so we may assume that we have an 
$n\times n$ binomial system. After dividing each binomial by a suitable 
monomial  we can then assume our system has the form $(x^{a_1}-c_1,\ldots,
x^{a_n}-c_n)$ for some $a_1,\ldots,a_n\!\in\!\Zn$ and 
$c_1,\ldots,c_n\!\in\!L^*$. Furthermore, via 
a monomial change of variables, we may in fact assume that 
$x^{a_i}\!=\!x^{d_i}_i$ for all $i$, for some choice of integers $d_1,
\ldots,d_n$. The latter reduction is routine, but we are 
unaware of a treatment in the literature allowing general fields. So we 
present a concise version below. 

For any integral matrix $A\!=\![a_{i,j}]\!\in\!\Z^{n\times n}$ with 
columns $a_1,\ldots,a_n$, let us write\linebreak 
$x^A\!=\!(x^{a_1},\ldots,x^{a_n})$ 
where the notation $x^{a_i}\!=\!x^{a_{1,i}}_1\cdots x^{a_{n,i}}_n$ is 
understood. It is easily checked that $x^{AB}\!=\!(x^A)^B$ for any $n\times n$ 
matrix $B$. 

Recall that an integral matrix 
$U\!\in\!\Z^{n\times n}$ is said to be {\em unimodular} if and only if its 
determinant is $\pm 1$. It is easily checked that the substitution $x\!=\!y^U$ 
induces an automorphism on $(L^*)^n$ that also preserves the number of roots 
with all coordinates having generalized phase $1$. From the classical theory of 
{\em Smith factorization} \cite{smith,storjophd}, one can always write 
$UAV\!=\!D$ for some unimodular $U$ and $V$, and a diagonal matrix $D$ with 
nonnegative diagonal entries $d_1,\ldots,d_n$.  

Applying the last two paragraphs to our binomial system $x^A-c$, 
we see that to count the maximal number of roots in $(L^*)^n$ (with all 
coordinates having generalized phase $1$) we may assume that 
our system is in fact $(x^{d_1}_1-c_1,\ldots,x^{d_n}_n-c_n)$. We thus 
obtain $Y_L(n,1)\!=\!Y_L(1,1)^n$ and, by Assertions (2), (1), (4), and (6) of 
Theorem \ref{thm:sum}, we are done. \qed 

\subsection{Proof of Theorem \ref{thm:hah}} 
The inequality $Y_\C(n,k)\!\geq\!Y_\R(n,k)$ is immediate 
since any real $(n+k)$-nomial $n\times n$ system is automatically 
a complex $(n+k)$-nomial $n\times n$ system. So we need only prove 
that $Y_\C(n,k)\!\leq\!Y_\R(n,k)$. To do the latter, 
it clearly suffices to show that for any $(n+k)$-nomial $n\times n$ 
system $G\!:=\!(g_1,\ldots,g_n)$ over $\C$, with $N$ non-degenerate roots 
in $\Rn_+$, we can find an $(n+k)$-nomial $n\times n$ system 
$F\!:=\!(f_1,\ldots,f_n)$ --- with all coefficients {\em real} --- having at 
least $N$ non-degenerate roots in $\Rn_+$. So, for all $i$, let us define 
$f_i\!:=\!e^{\sqrt{-1}t}g_i+e^{-\sqrt{-1}t}\bg_i$ where $\bar{(\cdot)}$ denotes 
complex conjugation, $\bg_i$ is the polynomial obtained from $g_i$ by 
conjugating all its coefficients, and $t\!\in\![0,2\pi)$ is a constant to be 
determined later. Clearly, for all $i$, 
the coefficients of $f_i$ are all real, and any exponent vector appearing 
in $f_i$ also appears in $g_i$. 

It is also clear that for any $\zeta\!\in\!\Rn_+$ with $G(\zeta)\!=\!0$ we 
have\\ 
\mbox{}\hfill  
$f_i(\zeta)\!=\!e^{\sqrt{-1}t}g_i(\zeta)+e^{-\sqrt{-1}t}\bg_i(\zeta)
\!=\!e^{\sqrt{-1}t}g_i(\zeta)+\overline{e^{\sqrt{-1}t}g_i(\zeta)}\!=\!0$. 
\hfill \mbox{}\\ 
So any root of $G$ in $\Rn_+$ is a root of $F$ in $\Rn_+$. 

Let $\jac(F)(\zeta)$ denote the Jacobian determinant of $F$ 
evaluated at $\zeta$, and assume now that $\zeta\!\in\!\Rn_+$ is a 
non-degenerate root of $G$. To see that $\zeta$ is also a non-degenerate 
root of $F$ (for a suitable choice of $t$), note that the multi-linearity of 
the determinant implies the following:\\  
\mbox{}\hfill 
$\jac(F)(\zeta)=\sum\limits_{s=(s_1,\ldots,s_n)\in\{\pm\}^n} 
e^{\sqrt{-1}(n_+(s)-n_-(s))t}\jac(g_{1,s_1},\ldots,g_{n,s_n})(\zeta)$,   
\hfill \mbox{}\\ 
where $n_\pm(s)$ is the number of $\pm$ signs in $s$, $g_{i,+}\!:=\!g_i$, 
and $g_{i,-}\!:=\!\bg_i$. In particular, we see that 
$\jac(F)(\zeta)\!=\!J\!\left(e^{\sqrt{-1}t}\right)$ for some $J\!\in\!\C\!
\left[x_1,\frac{1}{x_1}\right]$. Moreover, $J$ is not identically zero since 
the coefficient of $x^n_1$ is $\jac(G)(\zeta)\!\neq\!0$.
Clearly then, $J$ has at most $2n$ roots in $\Cs$ and thus there are at 
most $2n$ values of $t\!\in\![0,2\pi)$ for which $\jac(F)(\zeta)$ vanishes. 

Thus, assuming $G$ has $N$ non-degenerate roots in $\Rn_+$, 
$F$ fails to have at least $N$ non-degenerate roots in $\Rn_+$ for 
at most $2nN$ values of $t\!\in\![0,2\pi)$. \qed  

\subsection{Proof of Lemma \ref{lemma:k1}} 
Recall that in Section \ref{sec:proof} we wrote $R_n\!=\!A_n-B_n$ 
where $A_n$ and $B_n$ are suitable monomials. Assuming 
$n\!\geq\!3$ is odd we obtain the following:  \\  
\mbox{}\hfill $\displaystyle{A_n\!\left(\frac{16^{n-2}}{u}\right) =\frac{16^{n-2}}{u}
\prod_{i=1}^{\floor{n/2}}\left(1+4^{3-4i}\frac{4^{2n-4}}{u}\right)^2 
=\frac{16^{n-2}}{u}\prod_{i=1}^{\floor{n/2}}
\left(1+\frac{4^{2n-4i-1}}{u}\right)^2}$\hfill\mbox{}\\ 
\mbox{}\hfill $\displaystyle{=\frac{16^{n-2}}{u}\prod_{i=1}^{\floor{n/2}}
\left(\frac{4^{2n-4i-1}}{u}\left(1+4^{4i-2n+1}u\right)\right)^2 
=\frac{4^{2n-4}}{u}\cdot\frac{4^S}{u^{n-1}}\prod_{i=1}^{\floor{n/2}}
\left(1+4^{4i-2n+1}u\right)^2}$,\hfill\mbox{}\\ 
where $\displaystyle S\!=\!2\sum_{i=1}^{\floor{n/2}}(2n-4i-1)$.  
A minor calculation shows that $S+2n-4=(n-2)(n+1)$, so replacing 
$i$ by $\floor{n/2}-i+1$, we get\\ 
\mbox{}\hfill 
$\displaystyle A_n\!\left(\frac{16^{n-2}}{u}\right)=\left(\frac{4^{n-2}}{u}\right)^{n+1}u\prod_{i=1}^{\floor{n/2}}(1+4^{3-4i}u)^2=
\left(\frac{4^{n-2}}{u}\right)^{n+1}A_n(u)$. \hfill \mbox{}\\ 
An almost identical calculation proves the same transformation 
law for $B_n(u)$. Since\linebreak $R_n\!=\!A_n-B_n$, we thus 
obtain our transformation law for odd $n$. 

For even $n$, a similar calculation yields 
$A_n\!\left(\frac{16^{n-2}}{u}\right)
\!=\!\left(\frac{4^{n-2}}{u}\right)^{n+1}B_n(u)$ and\linebreak  
$B_n\!\left(\frac{16^{n-2}}{u}\right)
\!=\!\left(\frac{4^{n-2}}{u}\right)^{n+1}A_n(u)$. So  
we obtain $R_n\!\left(\frac{16^{n-2}}{u}\right)
\!=\!-\left(\frac{4^{n-2}}{u}\right)^{n+1}R_n(u)$ and thus 
the first assertion is proved.   

The final assertion follows immediately from our transformation law 
since $16^{n-2}/4^{n-2}\!=\!4^{n-2}$ and $(-4^{n-2}/4^{n-2})^{n+1}\!=\!-1$ 
for even $n$. \qed  

\subsection{Proof of Lemma \ref{lemma:k2}}  
To prove (a), merely observe that $R_n(0)\!=\!-\frac{1}{16}\!<\!0$ for all 
$n\geq 2$. 

To prove (b), the cases $n\!\leq\!4$ can be verified by direct 
computation. So let us assume $n\!\geq\!5$ and separate into two 
exclusive cases.  

\smallskip 
\noindent 
{\bf ($\ell$ even):} Let us first observe the following elementary inequality: 
\begin{eqnarray} 
\label{ineq:key} 
\prod^{(n-1)/2}_{i=1} \left(1-\frac{15/16}{1+256^{i-2}}\right)
& \geq & \frac{7}{200}\left(1+\frac{1}{4^{n-1}}\right) 
\text{ for all odd } n\!\geq\!3. 
\end{eqnarray} \\ 
Inequality (\ref{ineq:key}) follows easily by induction, after one 
first verifies the cases $n\!\in\!\{3,5,7\}$ directly. The identity 
$\frac{1+16z}{1+z}\!=\!16\left(1-\frac{15/16}{1+z}\right)$ then 
easily implies the following equality: 
\begin{eqnarray} 
\label{eq:key} 
\left(\frac{1+4^{2n-8}}{1+4^{2n-10}}\right)
\left(\frac{1+4^{2n-12}}{1+4^{2n-14}}\right)\cdots 
\left(\frac{1+4^{-2}}{1+4^{-4}}\right) = 16^{(n-1)/2}\prod^{(n-1)/2}_{i=1} 
\left(1-\frac{15/16}{1+256^{i-2}}\right) 
\end{eqnarray} 

Combining (\ref{ineq:key}) and (\ref{eq:key}) we then 
obtain, for any odd $n\!\geq\!5$:  
\begin{eqnarray*} 
\frac{A_n\!\left(4^{2n-7}\right)}{B_n\!\left(4^{2n-7}\right)}&=&\frac{4^{2n-7}\cdot 4^{2n-2}} 
{\left(\frac{1}{4}+4^{2n-7}\right)^2}\prod^{(n-1)/2}_{i=1} 
\left(1-\frac{15/16}{1+256^{i-2}}\right)^2  \\ 
& \geq & \frac{4^{2n-7}\cdot 4^{2n-2}} 
{\left(\frac{1}{4}+4^{2n-7}\right)^2}\frac{7^2}{200^2}
\left(1+\frac{1}{4^{n-1}}\right)^2
= \frac{4^{2n-7}\cdot 4^{2n-7}} 
{\left(\frac{1}{4}+4^{2n-7}\right)^2}\cdot \frac{4^5\cdot 7^2}{200^2}
\left(1+\frac{1}{4^{n-1}}\right)^2\\ 
& = & \left(\frac{1+\frac{1}{4^{n-1}}} 
{1+\frac{1}{4^{2n-6}}}\right)^2\cdot \frac{4^5\cdot 7^2}{200^2}\geq  
\frac{4^5\cdot 7^2}{200^2}=1.2544>1 
\end{eqnarray*}\\ 
We thus obtain  
\begin{eqnarray} 
\label{ineq:crux} 
A_\ell\!\left(4^{2\ell-7}\right) & > & B_\ell\!\left(4^{2\ell-7}\right)  
\text{for all odd } \ell\!\geq\!3 
\end{eqnarray} 

Recall that for any odd $n$,  
(i) $A_{n+1}(u)\!=\!A_n(u)\left(1+\frac{u}{4^{2n-1}}\right)^2$ and  
$B_{n+1}(u)\!=\!B_n(u)$, and (ii) 
$A_{n+2}(u)\!=\!A_n(u)\left(1+\frac{u}{4^{2n-1}}\right)^2$ and  
$B_{n+1}(u)\!=\!B_n(u)\left(1+\frac{u}{4^{2n+1}}\right)^2$. 
Combining the recurrences (i) and (ii) with Inequality (\ref{ineq:crux}), we 
then easily obtain by induction and re-indexing that 
$A_n(16^\ell/4)\!>\!B_n(16^\ell/4)$ for all $\ell\!\in\!\{0,\ldots, 
n-3\}$ with $\ell$ even. So we are done. \qed 

\medskip 
\noindent
{\bf ($\ell$ odd):} This case follows almost identically as the last case, 
save for minor 
changes in the indexing. In particular, one first uses Inequality 
(\ref{ineq:key}) to prove that 
$A_\ell\!\left(4^{2\ell-7}\right)\!<\!B_\ell\!\left(4^{2\ell-7}\right)$ for 
all {\em even} $\ell\!\geq\!4$. One then increases the subscript from 
$\ell$ to $n$ by induction, and re-indexes $\ell$, just as before. So we omit 
the details for brevity. \qed 

\section*{Acknowledgements} 
We thank Mart\'in Avenda\~{n}o, Qi Cheng, Pascal Koiran, Jeff Lagarias, Matt 
Papanikolas, Philippe P\'ebay, and Bjorn Poonen for useful discussions, some 
of them $p$-adic.  In particular, Example \ref{ex:slp} was inspired by a 
conversation between the second author and Bjorn Poonen at the Extensions of 
Hilbert's Tenth Problem workshop at the American Institute of Mathematics, 
March 21--25, 2005.  We also thank the anonymous referee for useful commentary. 
Special thanks go to Richard Lipton for pointing out his earlier result 
\cite{lipton} that cryptographic hardness of factoring implies a weakened 
version of the $\tau$-Conjecture. Finally, we thank the Fields Institute 
and the Laboratoire de l'Informatique du Parall\'elisme (at 
the \'Ecole Normale Sup\'erieure de Lyon) 
for their hospitality and great work environment while this paper was 
completed. 

We dedicate this paper in honor of Bernd 
Sturmfels' 50$\thth$ birthday. Happy 50 Bernd! 

\bibliographystyle{amsalpha}

\end{document}